\documentclass[10pt]{article}
\usepackage{amsfonts,color}
\usepackage{amssymb}
 \usepackage{footnote}
\usepackage{mathtools}
\usepackage{amsthm,wasysym}
\usepackage[english]{babel}
\usepackage{bbm}
\usepackage{colonequals}
\usepackage{cite}
\usepackage{setspace}

\usepackage{graphicx}
\usepackage{amsfonts}
\usepackage{hyperref}
\usepackage{subcaption}
\usepackage{amsmath}
\usepackage{nicefrac}
\usepackage{gensymb}
\usepackage{enumerate}
\usepackage[nameinlink,capitalize]{cleveref}

\usepackage[nottoc,numbib]{tocbibind}

\usepackage{pgf,tikz,pgfplots}
\pgfplotsset{compat=1.14}
\usepackage{mathrsfs}
\usetikzlibrary{arrows}

\DeclareMathOperator{\at}{\bigg\vert}
\newcommand{\vb}[1]{\mathbf{#1}}
\newcommand{\bm}[1]{\boldsymbol{#1}}

\DeclareMathOperator{\di}{\mathrm{div}}
\DeclareMathOperator{\dev}{\mathrm{dev}}
\DeclareMathOperator{\Di}{\mathrm{Div}}
\DeclareMathOperator{\D}{\mathrm{D}\hspace{-0.05cm}}
\DeclareMathOperator{\sym}{\mathrm{sym}}
\DeclareMathOperator{\spa}{\mathrm{span}}

\DeclareMathOperator{\skw}{\mathrm{skew}}
\DeclareMathOperator{\Le}{\mathit{L}^2}

\DeclareMathOperator{\Hone}{\mathit{H}^1}

\DeclareMathOperator{\Po}{\mathit{P}}

\DeclareMathOperator{\so}{\mathfrak{so}(3)}
\DeclareMathOperator{\Sph}{\mathbb{R}\cdot\bm{\mathbbm{1}}}
\DeclareMathOperator{\Sym}{\mathrm{Sym}}
\DeclareMathOperator{\GL}{\mathrm{GL}}

\DeclareMathOperator{\Ned}{\mathcal{N}}
\DeclareMathOperator{\Lag}{\mathcal{L}}
\DeclareMathOperator{\San}{\mathcal{T}}

\DeclareMathOperator{\one}{\bm{\mathbbm{1}}}
\newcommand{\Hd}[1]{\mathit{H}(\mathrm{div}{#1})}
\newcommand{\Hc}[1]{\mathit{H}(\mathrm{curl}{#1})}
\newcommand{\Hsc}[1]{\mathit{H}(\mathrm{sym}\, \mathrm{Curl}{#1})}

\newcommand{\HdD}[1]{\mathit{H}(\di\mathrm{Div}{#1})}
\newcommand{\jump}[1]{\ensuremath{[\![#1]\!]} }

\DeclareMathOperator{\Anti}{\mathrm{Anti}}
\DeclareMathOperator{\tr}{\mathrm{tr}}

\DeclareMathOperator{\curl}{\mathrm{curl}}
\DeclareMathOperator{\range}{\mathrm{range}}
\DeclareMathOperator{\Curl}{\mathrm{Curl}}

\DeclareMathOperator{\muma}{\mu_{\mathrm{macro}}}
\newcommand{\Lc}{L_\mathrm{c}}

\newcommand{\dd}{\mathrm{d}}
\DeclareMathOperator{\Pm}{\bm{P}}
\allowdisplaybreaks[1]

\makeindex

\newtheorem{remark}{Remark}
\usepackage{chngcntr}
\counterwithin{theorem}{section}
\counterwithin{lemma}{section}
\counterwithin{remark}{section}

   \makeatletter
\let\@fnsymbol\@arabic
\makeatother

\crefname{Problem}{Problem.}{Problem.}

\title{On $[\Hone]^{3 \times 3}$, $[\Hc{}]^3$ and $\Hsc{}$ finite elements for matrix-valued Curl problems}

\author{\normalsize{Adam Sky}\thanks{Corresponding author: Adam Sky, Institute of Structural Mechanics, Statics and Dynamics, Technische Universit\"at Dortmund, August-Schmidt-Str. 8, 44227 Dortmund, Germany, email: adam.sky@tu-dortmund.de}
	, \quad
	\normalsize{Ingo Muench}\thanks{Ingo Muench, Institute of Structural Mechanics, Statics and Dynamics, Technische Universit\"at Dortmund, August-Schmidt-Str. 8, 44227 Dortmund, Germany, email: ingo.muench@tu-dortmund.de}
	\quad
	and \quad
	\normalsize{Patrizio Neff}\thanks{Patrizio Neff, Chair for Nonlinear Analysis and Modelling, Faculty of Mathematics, Universit\"{a}t Duisburg-Essen,
		Thea-Leymann Str. 9, 45127 Essen, Germany, email: patrizio.neff@uni-due.de}
}

\begin{document}

\maketitle

\begin{abstract}
In this work we test the numerical behaviour of matrix-valued fields approximated by finite element subspaces of $[\Hone]^{3\times 3}$, $[\Hc{}]^3$ and $\Hsc{}$ for a linear abstract variational problem connected to the relaxed micromorphic model. The formulation of the corresponding finite elements is introduced, followed by numerical benchmarks and our conclusions.
The relaxed micromorphic continuum model reduces the continuity assumptions of the classical micromorphic model by replacing the full gradient of the microdistortion in the free energy functional with the Curl. This results in a larger solution space for the microdistortion, namely $[\Hc{}]^3$ in place of the classical $[\Hone]^{3 \times 3}$. The continuity conditions on the microdistortion can be further weakened by taking only the symmetric part of the Curl. As shown in recent works, the new appropriate space for the microdistortion is then $\Hsc{}$. The newly introduced space gives rise to a new differential complex for the relaxed micromorphic continuum theory. 
\\
\vspace*{0.25cm}
\\
{\bf{Key words:}}  relaxed micromorphic continuum, \and $\sym \Curl$ elements, \and N\'{e}d\'{e}lec elements, \and relaxed micromorphic complex, \and finite elements, \and approximation error. \\

\end{abstract}


\section{Introduction}
Today there exist various formulations of micromorphic continua \cite{Mindlin1968,Munch2011,Neff2014,Neff2009} and higher gradient theories \cite{highergradient} with the common goal of capturing micro-motions, that are unaccounted for in the classical Cauchy continuum theory.
The unifying characteristic of all micromorphic theories is the extension of the mathematical model with additional degrees of freedom (called the microdistortion) for the material point. Consequently, in the general micromorphic theory, as introduced by Eringen \cite{Eringen} and Mindlin \cite{Mindlin}, each material point is endowed with nine extra degrees of freedom given by the matrix $\Pm$, capturing affine displacements that are independent of the translational degrees of freedom, such as rotation or expansion of the micro-body. 

The relaxed micromorphic continuum \cite{Neff2014,Neff2019} differs from other micromorphic theories by reducing the continuity assumptions on these micro-motions. Instead of incorporating the full gradient $\D \Pm$ of the microdistortion into the free energy functional, the relaxed micromorphic continuum assumes only the Curl of the microdistortion to produce significant energies. As a result, the natural space for the microdistortion in the relaxed micromorphic continuum is $[\Hc{}]^3$. Further, the micro-dislocation, i.e. $\Curl \Pm$, remains a second-order tensor, in contrast to the general micromorphic theory where the full gradient yields a third order tensor $\D \Pm \in \mathbb{R}^{3\times 3 \times 3}$. Possible applications for the relaxed micromorphic theory, such as the simulation of metamaterials and bandgap materials, are demonstrated in \cite{Madeo2016, Madeo2018}. Further, closed-form solutions for specimen under shear, bending and torsion have already been derived in \cite{Rizzi2021,Rizzi2021b,rizzi2020analytical_2}. A first numerical implementation for a relaxed micromorphic model of antiplane shear can be found in \cite{Sky2021}, followed by an implementation for plain-strain models \cite{SSSN21} and lastly, an implementation of the full three-dimensional model \cite{Sky2}.

In recent works \cite{Starke2,Starke1, Lewintan2021,Lewintan2021c,Lewintan2021d,Lewintan2021b} it was shown, that the continuity assumptions on the microdistortion can be weakened further by considering solely the symmetric micro-dislocation $\sym \Curl \Pm$.
The corresponding space for the microdistortion is the Hilbert space $\Hsc{}$. A version of this space is also used for some formulations of the biharmonic equation \cite{Pauly2020}, with a restriction to trace-free tensors. 

In order to gain insight into the numerical behaviour of the microdistortion in these spaces and the potential meaning for the relaxed micromorphic continuum theory, we investigate three finite element formulations for $\bm{P} \in [\Hone]^{3\times 3}$, $\bm{P} \in [\Hc{}]^{3}$ and $\bm{P} \in \Hsc{}$ on an abstract problem derived from the relaxed micromorphic continuum. We compare standard Lagrangian, N\'ed\'elec- \cite{Nedelec1980,Ned2,Joachim2005,Zaglmayr2006,Sky2021} and the recently introduced $\Hsc{}$-base functions \cite{sander2021conforming}, respectively.

This paper is organized as follows: In the first two sections we introduce the relaxed micromorphic continuum and derive a related abstract problem. \cref{sec:fem} is devoted to the description of the finite element formulations. In \cref{seq:eg} we present numerical examples and investigate their behaviour. The last section presents our conclusions and outlook.

\section{A relaxed micromorphic complex}
\begin{figure}
    \centering
    \includegraphics[width=0.79\linewidth]{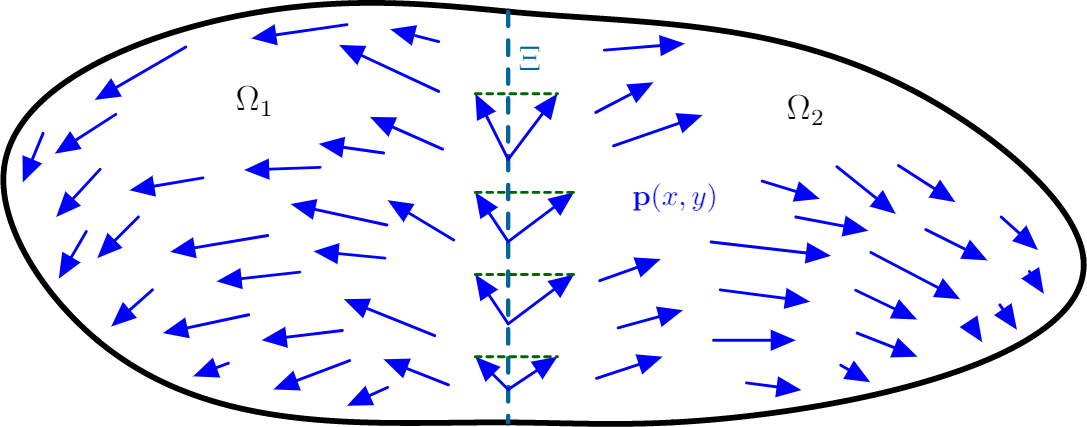}
    \caption{An example for a vector field $\vb{p} \in \Hc{,\Omega}$. Note that $\vb{p} \notin [\Hone(\Omega)]^2$ since on the interface $\Xi$, only the tangential component of the vector field is continuous.}
    \label{fig:hcurl_flux}
\end{figure}

The relaxed micromorphic continuum \cite{Neff2019,Neff2014} is described by its free energy functional, incorporating the gradient of the displacement field, the microdistortion and its Curl
\begin{gather}
	\begin{aligned} 
		I(\vb{u}, \bm{P}) =  \dfrac{1}{2} \int_{\Omega} & \langle \mathbb{C}_{\textrm{e}} \sym(\D\vb{u} - \bm{P}) , \, \sym(\D \vb{u} - \bm{P}) \rangle
		+  \langle \mathbb{C}_{\textrm{micro}} \sym\bm{P} , \, \sym\bm{P} \rangle \notag \\ 
		& + \langle \mathbb{C}_{\textrm{c}} \skw(\D \vb{u} - \bm{P}) , \, \skw (\D \vb{u} - \bm{P}) \rangle
		+ \mu_\text{macro} \, \Lc^2 \, \| \text{Curl}\bm{P}\|^2  \, \dd X  - \int_{\Omega} \langle \vb{u} , \, \vb{f} \rangle  - \langle \bm{P} , \, \bm{M} \rangle \, \dd X \, ,  
	\end{aligned}
	\label{eq:1} \\[2ex]
	\D \vb{u} = \begin{bmatrix}
		u_{1,1} & u_{1,2} & u_{1,3} \\
		u_{2,1} & u_{2,2} & u_{2,3} \\
		u_{3,1} & u_{3,2} & u_{3,3} 
	\end{bmatrix}
	\, , \quad
	\text{Curl}\bm{P} = \begin{bmatrix}
		\text{curl} \begin{bmatrix}
			P_{11} & P_{12} & P_{13} 
		\end{bmatrix} \\[1ex]
		\text{curl} \begin{bmatrix}
			P_{21} & P_{22} & P_{23} 
		\end{bmatrix} \\[1ex]
		\text{curl} \begin{bmatrix}
			P_{31} & P_{32} & P_{33} 
		\end{bmatrix}
	\end{bmatrix} \, , \quad 
	\text{curl}\vb{v} = \nabla \times \vb{v} \, , 
	\label{eq:full}
\end{gather}
where $\langle \cdot , \, \cdot \rangle$ denotes the scalar product on $\mathbb{R}^{3\times 3}$, $\Omega \subset \mathbb{R}^3$ is the domain and $\vb{u} : \Omega \subset \mathbb{R}^3 \to \mathbb{R}^3$ and $\bm{P}: \Omega \subset \mathbb{R}^3 \to \mathbb{R}^{3 \times 3}$ represent the displacement and the non-symmetric microdistortion, respectively. The volume forces and micro-moments are given by $\vb{f}$ and $\bm{M}$. Here, $\mathbb{C}_{\textrm{e}}$ and $\mathbb{C}_{\textrm{micro}}$ are standard fourth elasticity tensors and $\mathbb{C}_{\textrm{c}}$ is a positive semi-definite coupling tensor for rotations. The macroscopic shear modulus is denoted by $\mu_\text{macro}$ and the parameter $\Lc > 0$ represents the characteristic length scale motivated by the microstructure.
The corresponding Hilbert spaces and their respective traces read
\begin{align}
	\Hone(\Omega) &= \{ u \in \Le(\Omega) \; | \; \nabla u \in [\Le(\Omega)]^{3} \} \, , & \tr_{\Hone} u &= u \at_{\partial \Omega} \, , \\
	\Hc{,\Omega} &= \{ \vb{p} \in [\Le(\Omega)]^3 \; | \; \curl \vb{p} \in [\Le(\Omega)]^{3} \} \, , & \tr_{\Hc{}} \vb{p} &= \vb{p} \times \bm{\nu}  \at_{\partial \Omega} \, . 
\end{align}
Here, $\bm{\nu}$ denotes the unit outward normal on the surface. Existence and uniqueness of the relaxed micromorphic model using $\mathit{X} = [\Hone(\Omega)]^3 \times [\Hc{,\Omega}]^3$, where $[\Hc{,\Omega}]^3$ is to be understood as a row-wise matrix of the vectorial space, is derived in e.g. \cite{NeffExis,Neff_existence, GNMPR15}.
An example for a function belonging to $\Hc{}$ while not belonging to $\Hone$ is given in \cref{fig:hcurl_flux}.

Although the weak formulation of the relaxed micromorphic model does not represent a mixed formulation and therefore, does not require the use of the commuting de Rham complex for existence and uniqueness in the discrete case, it does introduce the so called \textbf{consistent coupling condition}\cite{dagostino2021consistent,Neff2019}
\begin{align}
	&\Pm \times \, \bm{\nu}  = \D \widetilde{\vb{u}} \times \bm{\nu}  \qquad \text{on} \quad  \Gamma_{\hspace{-0.5mm}D}  \, , 
\end{align}
where $\widetilde{\vb{u}}$ is the prescribed displacement field on the Dirichlet boundary $\Gamma_D = \Gamma_D^u = \Gamma_D^p$. The condition can only be satisfied exactly in the general discrete case, if commuting projections in the sense of a continuous-to-discrete de Rham complex are employed, compare with \cite{Demkowicz2000}. 
Further, when the characteristic length goes to infinity $\Lc \to \infty$, a mixed formulation introducing a new variable for the hyperstress $\bm{D} = \muma\,\Lc^2 \Curl \bm{P}$ is required to guarantee existence and uniqueness, as shown in \cite{Sky2021}.
Consequently, the appropriate complex for the relaxed micromorphic model is the classical de Rham complex \cite{Arnold2020ComplexesFC} in three dimensions \cref{fig:rham}, where 
\begin{align}
	\Hd{,\Omega} &= \{ \vb{p} \in [\Le(\Omega)]^3 \; | \; \di \vb{p} \in \Le(\Omega) \} \, , && \Di \Pm = \begin{bmatrix}
		\di\begin{bmatrix}
			P_{11} & P_{12} &  P_{13} 
		\end{bmatrix} \\[1ex]
	    \di\begin{bmatrix}
	    	P_{21} & P_{22} &  P_{23} 
	    \end{bmatrix} \\[1ex]
        \di\begin{bmatrix}
        	P_{31} & P_{32} &  P_{33} 
        \end{bmatrix} 
	\end{bmatrix} \, , && \di \vb{p} =  \nabla \cdot \vb{p} \, .
\end{align}
\begin{figure}
	\centering
	\begin{tikzpicture}[scale = 0.6][line cap=round,line join=round,>=triangle 45,x=1.0cm,y=1.0cm]
		\clip(3,7) rectangle (24.5,10);
		\draw (3,9) node[anchor=north west] {$[\mathit{H}^1(\Omega)]^3$};
		\draw (8.5,9) node[anchor=north west] {$[\Hc{, \Omega}]^3$};
		\draw [->,line width=1.5pt] (5.7,8.5) -- (8.5,8.5);
		\draw (6.5,9.5) node[anchor=north west] {$\D$};
		\draw [->,line width=1.5pt] (12.2,8.5) -- (15.2,8.5);
		\draw (12.8,9.5) node[anchor=north west] {$\Curl$};
		\draw (15.3,9) node[anchor=north west]  {$[\Hd{,\Omega}]^3$};
		\draw [->,line width=1.5pt] (18.8,8.5) -- (21.8,8.5);
		\draw (19.4,9.5) node[anchor=north west] {$\Di$};
		\draw (21.8,9) node[anchor=north west] {$[\Le(\Omega)]^3$};
	\end{tikzpicture}
    \caption{The classical de Rham exact sequence.
    The range of each operator is exactly the kernel of the next operator in the sequence, assuming a contractible domain $\Omega$.}
    \label{fig:rham}
\end{figure} 

The continuity assumptions on the microdistortion field $\Pm$ can be further reduced by taking only its symmetric part $\sym\Curl$ instead of the full Curl \cite{Lewintan2021,Lewintan2021b,Lewintan2021c, Lewintan2021d,Starke2}
\begin{align} 
	I(\vb{u}, \bm{P}) =  \, \dfrac{1}{2} \int_{\Omega} & \langle \mathbb{C}_{\textrm{e}} \sym(\D \vb{u} - \bm{P}) , \, \sym(\D \vb{u} - \bm{P}) \rangle
	+  \langle \mathbb{C}_{\textrm{micro}} \sym\bm{P} , \, \sym\bm{P} \rangle
	\label[Problem]{eq:sym}  \\ 
	& + \langle \mathbb{C}_{\textrm{c}} \skw(\D \vb{u} - \bm{P}) , \, \skw (\D \vb{u} - \bm{P}) \rangle
	+ \mu_\text{macro} \, \Lc^2 \, \| \sym \text{Curl}\bm{P} \|^2 \, \dd X  - \int_{\Omega} \langle \vb{u} , \, \vb{f} \rangle  - \langle \bm{P} , \, \bm{M} \rangle \, \dd X \, . \notag
\end{align}
Since $\| \sym \Curl \Pm \|^2 \leq \| \Curl \Pm \|^2$, this is a considerably weaker formulation than \cref{eq:full}.

\begin{figure}
    \centering
    \includegraphics[width=0.59\linewidth]{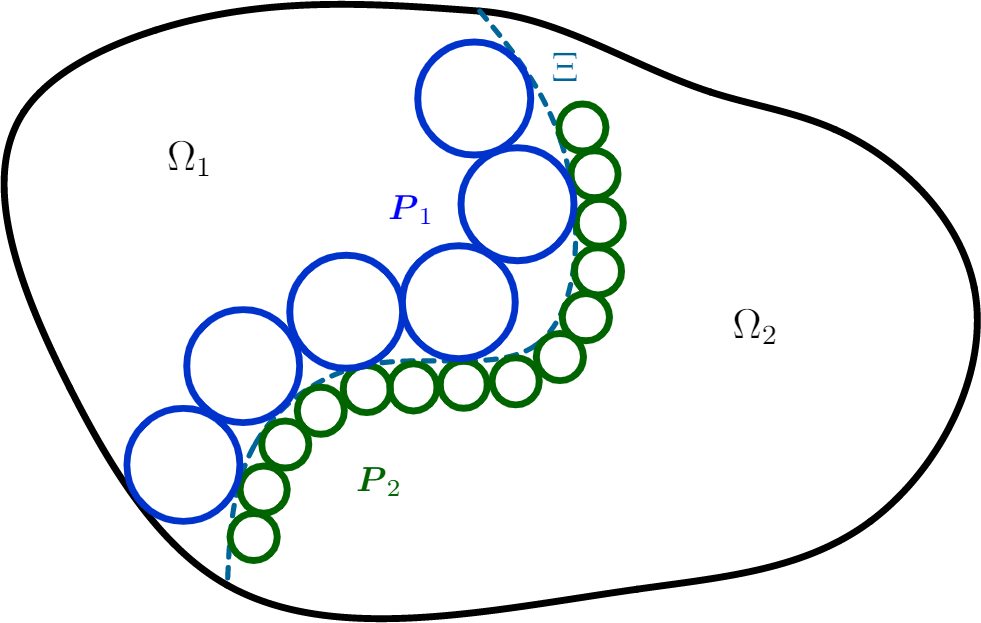}
    \caption{Depiction of a field $\Pm = \Pm_1 \cup \Pm_2 \in \Hsc{, \Omega}{}$, such that $\Pm \notin [\Hc{}]^3$. The circles illustrate the intensity of the spherical part of $\bm{P}$, which is discontinuous along the dashed line. } 
    \label{fig:hsymcurl_field}
\end{figure}

The appropriate Hilbert space for $\Pm$ now reads
\begin{align}
	\Hsc{,\Omega} &= \{ \Pm \in [\Le(\Omega)]^{3 \times 3} \; | \; \sym \Curl \Pm \in [\Le(\Omega)]^{3\times 3} \} \, , && \tr_{\Hsc{}} \Pm = \sym[\Pm \, \Anti(\bm{\nu})^T] \at_{\partial \Omega} \, ,
\end{align}
where $\Anti(\cdot)$ generates the skew symmetric matrix from a vector $\bm{\nu} \in \mathbb{R}^3$
\begin{align}
	& \Anti(\bm{\nu}) = \begin{bmatrix}
		0 & -\nu_3 & \nu_2 \\
		\nu_3 & 0  & -\nu_1 \\
		-\nu_2 & \nu_1  & 0 
	\end{bmatrix} \in \so  \, , && \Anti(\bm{\nu}) \vb{v} = \bm{\nu} \times \vb{v} \, , \quad \vb{v} \in \mathbb{R}^3 \, ,
\end{align}
The existence of minimizers for \cref{eq:sym} is shown in \cite{Lewintan2021,Lewintan2021b}.
The difference in smoothness between $[\Hc{}]^3$ and $\Hsc{}$ can be observed when considering spherical matrix fields $\Pm = p \, \one \in \Sph$
\begin{align}
	&\tr_{\Hc{}} \bm{P} = \Pm \, \Anti(\bm{\nu})^T \at_{\partial \Omega} = p \one \, \Anti(\bm{\nu})^T \at_{\partial \Omega} = p \Anti(\bm{\nu})^T \at_{\partial \Omega} & \forall \bm{P} \in \Sph \, , \\
	&\tr_{\Hsc{}} \bm{P} = \sym[\Pm \, \Anti(\bm{\nu})^T] \at_{\partial \Omega} = \sym[p \Anti(\bm{\nu})^T] \at_{\partial \Omega} = 0 & \forall \bm{P} \in \Sph \, .
	\label{eq:symcurltr}
\end{align}
The latter identity is evident due to $p \Anti(\bm{\nu}) \in \so$ and $\ker(\sym) = \so$. Consequently, the $\Hsc{}$ space captures discontinuous spherical tensor fields. A possible interpretation for the kinematics of such a field is depictable under the assumption that the field $\Pm$ represents a micro-strain field in the domain. In which case, material points can undergo discontinuous dilatation, see \cref{fig:hsymcurl_field}.

This new formulation gives rise to a corresponding complex, designated here the relaxed micromorphic complex, see \cref{fig:complex}, where the $\HdD{}$ space is defined as
\begin{align}
	\HdD{,\Omega} &= \{ \Pm \in [\Le(\Omega)]^{3 \times 3} \; | \; \di \Di \Pm \in \Le(\Omega) \} \, .
\end{align}
The complex can be seen as an extension of the $\di \Di$-sequence \cite{Pauly2020,Hu2021}, in which only trace-free (deviatoric) gradient fields $\dev \D \vb{u}$ are concerned , where $\dev \bm{X} = \bm{X} - \dfrac{1}{3} \tr (\bm{X}) \one$.
In fact, the right half of the complex is the same due to 
\begin{align}
	\range(\sym \Curl \dev) = \range(\sym \Curl) \subset \ker(\di \Di) \, ,
\end{align}  
since (for a full derivation see \cref{ap:A})
\begin{align}
		&\di \Di \sym \Curl \Pm = 0 \, , 
\end{align}
and 
\begin{align}
	\sym \Curl \Pm = \sym \Curl \left[ \dev \Pm + \dfrac{1}{3} \tr(\Pm) \, \one \right] = \sym \Curl \dev \Pm  \qquad \forall \Pm \in \mathbb{R}^{3 \times 3} \, .
\end{align}
The left side of the complex \cite{Lewintan2021} is derived from (see \cref{ap:A})
\begin{align}
	\ker(\sym \Curl) = \D \, [\Hone]^3 \cup \Sph \, . 
	\label{eq:kersymcurl}
\end{align}
As shown in \cite{Pauly2020}, the $\di \Di$-complex is an exact sequence for a topologically trivial domain. Since the right half of the relaxed micromorphic complex is the same as the right half of the $\di \Di$-complex and due to the exactness of \cref{eq:kersymcurl} as derived in \cite{Lewintan2021}, the relaxed micromorphic complex is also an exact sequence for a topologically trivial domain.
\begin{remark}
    The relaxed micromorphic complex serves to better understand the behaviour of the model with respect to the characteristic length $\Lc$. Let $\Lc \to \infty$, the next space in the sequence $\HdD{}$ is needed to approximate the hyperstress field
    \begin{align*}
        \bm{D} = \muma\,\Lc^2 \sym \Curl \bm{P} \in \HdD{, \Omega} \cap \Sym(3)
        \, , 
    \end{align*}
    and to express a stable mixed formulation, compare with \cite{Sky2021}. The construction of $\HdD{}$-conforming finite elements is beyond the scope of this work.
\end{remark}

\begin{figure}
	\centering
	\begin{tikzpicture}[scale = 0.6][line cap=round,line join=round,>=triangle 45,x=1.0cm,y=1.0cm]
		\clip(3,5) rectangle (28,10);
		\draw (3,9) node[anchor=north west] {$[\mathit{H}^1(\Omega)]^3$};
		\draw (8.5,9) node[anchor=north west] {$\Hsc{,\Omega}$};
		\draw [->,line width=1.5pt] (5.7,8.5) -- (8.5,8.5);
		\draw (6.5,9.5) node[anchor=north west] {$\D$};
		\draw (4.,6.4) node[anchor=north west] {$\Sph$};
		\draw [line width=1.5pt] (5.7,6) -- (10,6);
		\draw [->,line width=1.5pt] (10,6) -- (10,8);
		\draw (6.5,7) node[anchor=north west] {$\text{id}$};
		\draw (6.5,8) node[anchor=north west] {$\cup$};
		\draw [->,line width=1.5pt] (13,8.5) -- (16,8.5);
		\draw (13,9.5) node[anchor=north west] {$\sym \Curl$};
		\draw (16,9) node[anchor=north west] {$\HdD{,\Omega} \cap \sym(3)$};
		\draw [->,line width=1.5pt] (22.5,8.5) -- (25.5,8.5);
		\draw (22.7,9.5) node[anchor=north west] {$\di \Di$};
		\draw (25.5,9) node[anchor=north west] {$\Le(\Omega)$};
	\end{tikzpicture}
	\caption{The relaxed micromorphic complex for the microdistortion $\Pm$ and the corresponding curvature $\sym\Curl\Pm$. The kernel $\ker(\Hsc{,\Omega})$ is given fully by gradients of $\D \, [\Hone(\Omega)]^3$ and spherical tensors $\mathbb{R} \cdot \one$. The kernel $\ker(\HdD{,\Omega})$ is given by $\sym \Curl \Hsc{, \Omega}$. Finally, $\di \Di \HdD{, \Omega}$ is a surjection onto $\Le(\Omega)$. For all these relations we assume a contractible domain.}
	\label{fig:complex}
\end{figure} 

\section{Abstract variational problem}
In order to compare the behaviour of linear finite elements on $[\Hone]^{3 \times 3}$, $[\Hc{}]^3$ and $\Hsc{}$, the following abstract variational problem is introduced, drawing characteristics from the relaxed micromorphic model  
\begin{align}
		I(\bm{P}) &= \dfrac{1}{2} \int_\Omega \| \sym \bm{P} \|^2 + \| \sym \Curl \Pm \|^2  \, \dd X - \int_\Omega \langle \bm{P} ,\, \bm{M} \rangle \, \dd X \quad \to \quad \text{min} \, , \label{eq:abs} \\
		\delta I &= \int_\Omega \langle \sym \bm{P} ,\, \sym \delta \bm{P} \rangle + \langle \sym \Curl \bm{P} ,\, \sym \Curl (\delta\bm{P}) \rangle - \langle \delta \bm{P} ,\, \bm{M} \rangle \, \dd X = 0 \, .
		\label{eq:weak}
\end{align}
The associated strong form is derived by partial integration 
\begin{align}
	\sym \bm{P} + \Curl(\sym \Curl \bm{P}) &= \bm{M} && \text{in} \quad \Omega \, , \label{eq:coupled} \\
	\sym[\Pm \, \Anti(\bm{\nu})^T] &= \sym[\widetilde{\Pm} \, \Anti(\bm{\nu})^T] && \text{on} \quad \Gamma_{\hspace{-0.5mm}D} \, , \\
	\sym \Curl \bm{P} \times \bm{\nu} &= 0  && \text{on} \quad \Gamma_{\hspace{-0.5mm}N} \, ,
\end{align}
where $\widetilde{\bm{P}}$ is the prescribed field on the Dirichlet boundary $\Gamma_D$ and $\Gamma_N$ represents the Neumann boundary.
\begin{remark}
    Note that the variational problem resembles a vectorial Maxwell \cite{Joachim2005} problem 
    $\vb{v} + \curl\curl\vb{v} = \vb{m}$ for some vector field $\vb{v}: \Omega \mapsto \mathbb{R}^3$. However, in \cref{eq:coupled} the vectorial rows of $\Pm$ are highly coupled and the symmetrization of the Curl operator intervenes.
\end{remark}
For a predefined field $\widetilde{\bm{P}}$ where the entire boundary is prescribed ($\partial \Omega = \Gamma_{\hspace{-0.5mm}D}$) the micro-moment $\bm{M}$ is given by
\begin{align}
	\bm{M} = \sym \widetilde{\bm{P}} + \Curl(\sym \Curl \widetilde{\bm{P}}) \, , 
	\label{eq:microm}
\end{align}
and the analytical solution is $\bm{P} = \widetilde{\bm{P}}$.
Due to the generalized Korn's inequality \cite{Lewintan2021, Neff2012,Lewintan2021b,Lewintan2021c,Lewintan2021d,NEFF20151267}, the variational problem \cref{eq:abs} and the weak form \cref{eq:weak} are well-posed in the space $\Hsc{, \Omega}$.

\section{Finite element formulations}\label{sec:fem}
\begin{figure}
	\centering
	\begin{tikzpicture}
		\begin{axis}
			[
			width=30cm,height=17cm,
			view={50}{15},
			enlargelimits=true,
			xmin=-1,xmax=2,
			ymin=-1,ymax=2,
			zmin=-1,zmax=2,
			domain=-10:10,
			axis equal,
			hide axis
			]
			\draw[->, line width=1.pt, color=black](0., 0., 0.)--(1.5,0.,0.);
			\draw[color=black] (1.6,0,0) node[] {$\xi$};
			\draw[->, line width=1.pt, color=black](0., 0., 0.)--(0.,1.5,0.);
			\draw[color=black] (0.,1.6,0) node[] {$\eta$};
			\draw[->, line width=1.pt, color=black](0., 0., 0.)--(0.,0.,1.5);
			\draw[color=black] (0.,0.,1.6) node[] {$\zeta$};
			\addplot3[color=teal][
			line width=1.pt,
			mark=*
			]
			coordinates {
				(0,0,0)(1,0,0)(0,1,0)(0,0,0)
			};
			\addplot3[color=teal][
			line width=1.pt,
			mark=*
			]
			coordinates {
				(0,0,0)(0,0,1)
			};
			\addplot3[color=teal][line width=1.pt]coordinates {(1,0,0)(0,0,1)};
			\addplot3[color=teal][line width=1.pt]coordinates {(0,1,0)(0,0,1)};
			\draw[color=teal] (0,0,0) node[anchor=south east] {$_{(1)}$};
			\draw[color=teal] (1,0,0) node[anchor=north east] {$_{(2)}$};
			\draw[color=teal] (0,1,0) node[anchor=north west] {$_{(3)}$};
			\draw[color=teal] (0,0,1) node[anchor=north east] {$_{(4)}$};
			\fill[opacity=0.1, teal] (axis cs: 0,0,0) -- (axis cs: 1,0,0) -- (axis cs: 0,1,0) -- (axis cs: 0,0,1) -- cycle;
			
			\addplot3[color=teal][
			line width=1.pt,
			mark=*
			]
			coordinates {
				(2.2,2,0.5)(2.5,3,0.5)(3,2,0.2)(2,3,1.)
			};
		\addplot3[color=teal][line width=1.pt]coordinates {(2.2,2,0.5)(2,3,1.)};
		\addplot3[color=teal][line width=1.pt]coordinates {(2.5,3,0.5)(2,3,1.)};
		\addplot3[color=teal][line width=1.pt]coordinates {(2.2,2,0.5)(3,2,0.2)};
		\fill[opacity=0.1, teal] (axis cs: 2.2,2,0.5) -- (axis cs: 3,2,0.2) -- (axis cs: 2.5,3,0.5) -- (axis cs: 2,3,1.) -- cycle;
		\draw[color=teal] (axis cs: 2.2,2,0.5) node[anchor=east] {$\vb{x}_3$};
		\draw[color=teal] (axis cs: 3,2,0.2) node[anchor=east] {$\vb{x}_1$};
		\draw[color=teal] (axis cs: 2.5,3,0.5) node[anchor=west] {$\vb{x}_2$};
		\draw[color=teal] (axis cs: 2,3,1.) node[anchor=east] {$\vb{x}_4$};
		\draw[->, line width=1.pt, color=black](1.5, 1.5, 0.)--(2.5, 1.5,0.);
		\draw[color=black] (2.6, 1.5,0.) node[] {$x$};
		\draw[->, line width=1.pt, color=black](1.5, 1.5, 0.)--(1.5, 2.5,0.);
		\draw[color=black] (1.5, 2.6,0.) node[] {$y$};
		\draw[->, line width=1.pt, color=black](1.5, 1.5, 0.)--(1.5, 1.5,1.);
		\draw[color=black] (1.5, 1.5,1.1) node[] {$z$};
		
		\addplot3[smooth,color=blue][->, line width=1.pt]coordinates {(0.7,0.7,1.2)(1.2,1.2,1.4)(1.7,1.7,1.2)};
		\draw[color=blue] (1.2,1.2,1.4) node[anchor=south] {$\vb{x}(\xi, \, \eta , \, \zeta)$};
		\end{axis}
	\end{tikzpicture}
    \caption{Affine mapping from the reference element to the physical domain.}
    \label{fig:map} 
\end{figure}
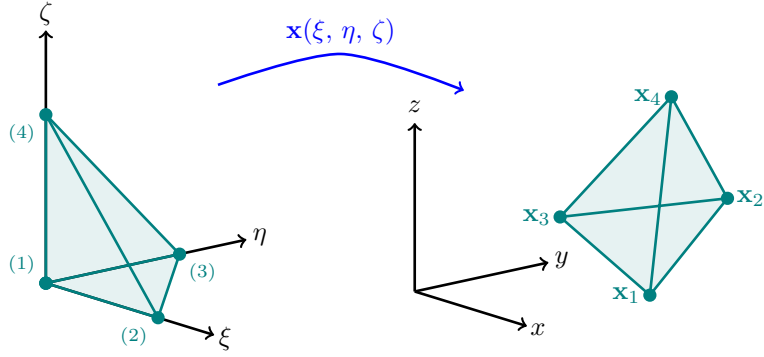

In the following we introduce finite elements using Voigt notation.
As a result, the microdistortion and micro-moment fields $\Pm$ and $\bm{M}$ are given by nine-dimensional vectors and the symmetry operator by the nine-by-nine matrix $\mathbb{S}$
\begin{align}
	& \Pm = \begin{bmatrix}
		P_{11} & P_{12} & P_{13} & P_{21} & P_{22} & P_{23} & P_{31} & P_{32} & P_{33}
	\end{bmatrix}^T \, , &&\bm{\Pm}, \, \bm{M} \in \mathbb{R}^9 \, , && \mathbb{S} \in  \mathbb{R}^{9 \times 9} \, .
\end{align}
All formulations apply to a linear tetrahedral element with a barycentric mapping (see \cref{fig:map})
\begin{align}
	&\vb{x}(\xi , \, \eta , \, \zeta) = (1 - \xi - \eta - \zeta) \, \vb{x}_1 + \xi \, \vb{x}_1 + \eta \, \vb{x}_2 +  \zeta \, \vb{x}_3 \, , &&
	\bm{J} = \D \vb{x} = \begin{bmatrix}
		\vb{x}_2 - \vb{x}_1 && \vb{x}_3 - \vb{x}_1 && \vb{x}_4 - \vb{x}_1 
	\end{bmatrix} \, ,
\end{align}
where $\bm{J}$ is the Jacobi matrix.

\subsection{Lagrangian $[\Hone]^{3\times 3}$-element}
The functionals of Lagrangian finite elements are defined by point-wise evaluation
\begin{align}
	l_i (u) \at_{x_j} = \delta_{ij} \, u \at_{x_j} \, ,
\end{align}
where $\delta_{ij}$ is the Kronecker delta. 
For the construction of the linear Lagrangian element, the evaluation of the degrees of freedom for the polynomial space $\Po^1 = \spa\{1, \, \xi , \, \eta , \, \zeta \}$ on the reference element yields the same barycentric base functions used for the element mapping 
\begin{align}
	&n_1(\xi , \, \eta , \, \zeta) = 1 - \xi - \eta - \zeta \, , && n_2(\xi , \, \eta , \, \zeta) = \xi \, , \\
	&n_3(\xi , \, \eta , \, \zeta) = \eta \, , && n_4(\xi , \, \eta , \, \zeta) = \zeta \, .
\end{align}
The entire ansatz matrix for the interpolation of a nine-dimensional vector can be built accordingly
\begin{align}
	&\bm{N} = \begin{bmatrix}
		n_1 \one_9 & n_2 \one_9 & n_3 \one_9 & n_4 \one_9
	\end{bmatrix} \in \mathbb{R}^{9 \times 36} \, , && \one_9 \in \mathbb{R}^{9\times 9} \, , && [\one_9]_{ij} = \delta_{ij} \, , && i,j \in \{1,2,\dots,9\} \, ,
\end{align}
where $\one_9$ is the nine-dimensional unit matrix. 
The curl of a row of $\Pm$ with respect to the physical domain is retrieved following the standard chain rule
\begin{align}
	\curl_x \vb{p}_i = \bm{J}^{-T} \curl_{\xi}  \vb{p}_i = \bm{J}^{-T} \Anti(\nabla_{\xi}) \, \vb{p}_i \, .
\end{align}
The matrix Curl operator is defined by the row-wise application of the curl operator.
The resulting element stiffness matrix reads
\begin{align}
	\bm{K}_\text{lag} = \int_\Omega (\bm{N}^T\mathbb{S}\, \bm{N} + \Curl (\bm{N})^T \, \mathbb{S} \, \Curl \bm{N}) \, \dd X \in \mathbb{R}^{36 \times 36} \, ,
\end{align}
where $\mathbb{S}$ defines the equivalent matrix symmetry operator in Voigt notation.
The Lagrangian element belongs to $\Lag^1 \subset \Hone{}$ and its interpolant produces the error estimate 
\begin{align}
    \| u - \Pi_{\Lag} u \|_{\Le} \preceq h^2 \,  |u|_{\mathit{H}^2} \, ,
\end{align}
where $u$ is the exact solution, $\Pi_{\Lag}$ is the Lagrangian interpolant and $h$ is a measure of the element's size.

\subsection{N\'ed\'elec $[\Hc{}]^{3}$-element}
\begin{figure}
	\centering
	\begin{subfigure}{0.15\linewidth}
		\includegraphics[width=1.0\linewidth]{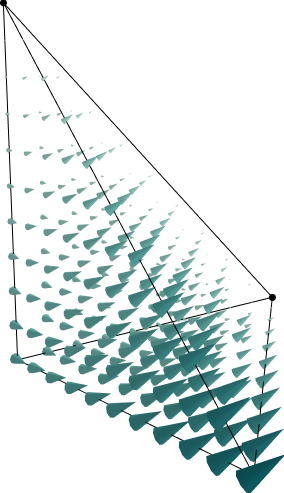}
		\caption{$\bm{\vartheta}_1$} 
	\end{subfigure}
	\begin{subfigure}{0.15\linewidth}
		\includegraphics[width=1.0\linewidth]{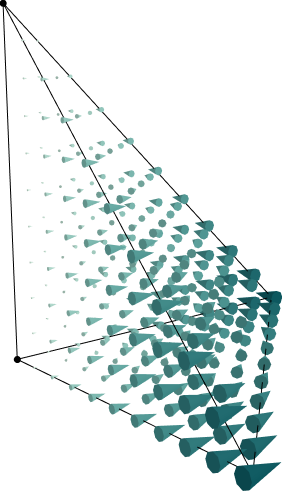}
		\caption{$\bm{\vartheta}_2$}
	\end{subfigure}
	\begin{subfigure}{0.15\linewidth}
		\includegraphics[width=1.0\linewidth]{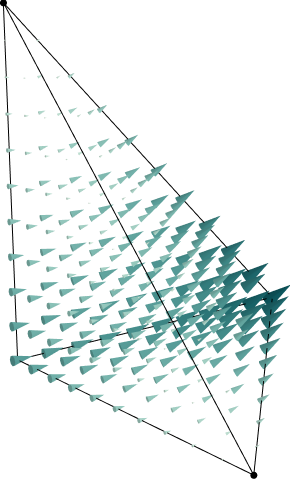}
		\caption{$\bm{\vartheta}_3$}
	\end{subfigure}
	\begin{subfigure}{0.15\linewidth}
		\includegraphics[width=1.0\linewidth]{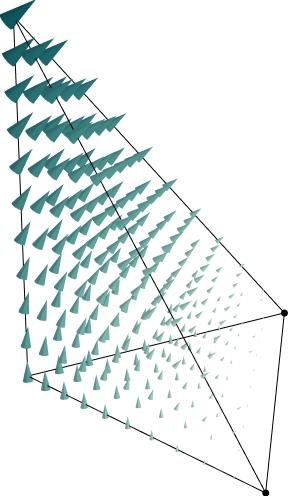}
		\caption{$\bm{\vartheta}_4$}
	\end{subfigure}
	\begin{subfigure}{0.15\linewidth}
		\includegraphics[width=1.0\linewidth]{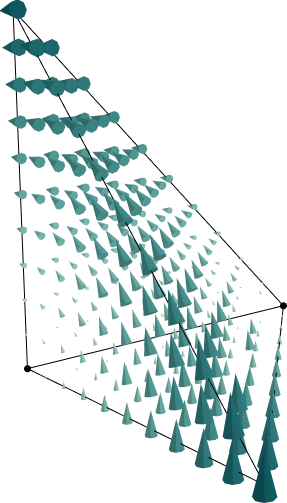}
		\caption{$\bm{\vartheta}_5$}
	\end{subfigure}
    \begin{subfigure}{0.15\linewidth}
    	\includegraphics[width=1.0\linewidth]{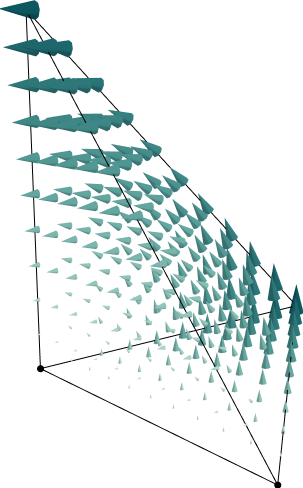}
    	\caption{$\bm{\vartheta}_6$}
    \end{subfigure}
	\caption{Linear N\'ed\'edelec base functions on the reference tetrahedron.}
	\label{fig:ned}
\end{figure}

N\'ed\'elec \cite{Nedelec1980} finite elements are conforming in $\Hc{}$. This is achieved by controlling the tangential projections of the base functions on the edges of a finite element. The lowest degrees of freedom, namely edge type, are defined as integrals along the curve of the element's edge
\begin{align}
	 l_{ij}(\vb{p}) = \int_{s_i} q_j \langle \vb{p} , \, \bm{\tau} \rangle \, \dd s \qquad \forall \, q_j \in \Po^{p-1}(s_i) \, .
\end{align}
Here $q_j$ are test functions. The adequate polynomial space is defined in \cite{Nedelec1980} as 
\begin{align}
	&\mathit{R}^p = [\Po^{p-1}]^3 \oplus \mathit{S}^p \, , && \mathit{S}^p = \{ \vb{p} \in [\widetilde{\Po}^p]^3 \; | \; \langle \vb{p} , \, \vb{x} \rangle = 0 \} \, ,
\end{align} 
where $\widetilde{\Po}$ is the space of homogeneous polynomials.
For linear polynomials $p = 1$ the space reads
\begin{align}
	&\mathit{R}^1 = \spa \left\{ \begin{bmatrix}
		1 \\ 0 \\ 0 
	\end{bmatrix} , \, \begin{bmatrix}
  0 \\ 1 \\ 0 
\end{bmatrix} \, ,
\begin{bmatrix}
	0 \\ 0 \\ 1 
\end{bmatrix} \, ,
\begin{bmatrix}
	0 \\ \zeta \\ -\eta 
\end{bmatrix} \, ,
\begin{bmatrix}
	-\zeta \\ 0 \\ \xi 
\end{bmatrix} \, ,
\begin{bmatrix}
	\eta \\ -\xi \\ 0 
\end{bmatrix} 
\right\} \, , && \dim \mathit{R}^1 = 6 \, .  
\end{align} 
The latter describes the polynomial space for linear N\'ed\'elec elements of the first type $\mathcal{N}_I^0$. Evaluating the degrees of freedom on the reference tetrahedron using the test functions $q_j = 1$ yields the corresponding base functions (see \cref{fig:ned})
\begin{align}
	&\bm{\vartheta}_1 = \begin{bmatrix}
		1 - \eta - \zeta \\
		\xi \\
		\xi
	\end{bmatrix} \, , &&
	\bm{\vartheta}_2 = \begin{bmatrix}
		- \eta \\
		\xi \\
		0
	\end{bmatrix} \, , &&
	\bm{\vartheta}_3 = \begin{bmatrix}
		\eta \\
		1-\xi-\zeta \\
		\eta
	\end{bmatrix} \, , \notag \\
	&\bm{\vartheta}_4 = \begin{bmatrix}
		\zeta \\
		\zeta \\
		1 - \xi - \eta
	\end{bmatrix} \, , &&
	\bm{\vartheta}_5 = \begin{bmatrix}
		- \zeta \\
		0 \\
		\xi
	\end{bmatrix} \, , &&
	\bm{\vartheta}_6 = \begin{bmatrix}
		0 \\
		-\zeta \\
		\eta
	\end{bmatrix} \, . 
\end{align}
\begin{remark}
	The zero power notation of the space $\mathcal{N}_I^0$ reminds that the base functions generate constant tangential projections on the edge. An alternative construction using the full polynomial space is given by the N\'ed\'elec elements of the second type $\mathcal{N}_{II}^1$ \cite{Ned2}.
\end{remark}
The base functions are vectors defined on the reference element. In order to map them to the physical domain we employ the covariant Piola transformation
\begin{align}
	\bm{\theta}_i = \bm{J}^{-T} \, \bm{\vartheta}_i \, .
\end{align}
Further, the curl of vectors undergoing a covariant Piola mapping is given by the contravariant Piola transformation
\begin{align}
	\curl_x \bm{\theta}_i = \dfrac{1}{\det \bm{J}} \bm{J} \curl_\xi \bm{\vartheta}_i \, . 
	\label{eq:piola}
\end{align}
\begin{remark}
	Note that Piola transformations guarantee a consistent projection on the element's boundaries in terms of size. However, the transformation does not control whether the tangential projections of neighbouring elements are parallel or anti-parallel, since the vectorial base functions are mapped to the physical domain separately by the corresponding Jacobi matrix of each element. Consequently, a correction function is employed to assert consistency, compare with \cite{Sky2021}.
\end{remark}

For the construction of the finite element we define a corresponding ansatz matrix
\begin{align}
	\bm{\Theta} = \begin{bmatrix}
		\bm{\theta}_1 & \vb{o} & \vb{o} & \bm{\theta}_2 & \vb{o} & \vb{o} & \bm{\theta}_3 & \vb{o} & \vb{o} & \bm{\theta}_4 & \vb{o} & \vb{o} & \bm{\theta}_5 & \vb{o} & \vb{o} & \bm{\theta}_6 & \vb{o} & \vb{o} \\
		\vb{o} & \bm{\theta}_1 & \vb{o} & \vb{o} & \bm{\theta}_2 & \vb{o} & \vb{o} & \bm{\theta}_3 & \vb{o} & \vb{o} & \bm{\theta}_4 & \vb{o} & \vb{o} & \bm{\theta}_5 & \vb{o} & \vb{o} & \bm{\theta}_6 & \vb{o}  \\
		\vb{o} & \vb{o} & \bm{\theta}_1 & \vb{o} & \vb{o} & \bm{\theta}_2 & \vb{o} & \vb{o} & \bm{\theta}_3 & \vb{o} & \vb{o} & \bm{\theta}_4 & \vb{o} & \vb{o} & \bm{\theta}_5 & \vb{o} & \vb{o} & \bm{\theta}_6
	\end{bmatrix} \in \mathbb{R}^{9 \times 18} \, ,
\end{align}
where $\vb{o} = \begin{bmatrix}
	0 & 0 & 0
\end{bmatrix}^T$ is a three-dimensional vector of zeros.
The Curl of the microdistortion $\Pm$ is calculated using \cref{eq:piola} for each base function
\begin{align}
	\Curl \bm{\Theta} = \begin{bmatrix}
		\curl \bm{\theta}_1 & \vb{o} & \vb{o} &  & \curl\bm{\theta}_6 & \vb{o} & \vb{o} \\
		\vb{o} & \curl \bm{\theta}_1 & \vb{o}  & \cdots & \vb{o} & \curl\bm{\theta}_6 & \vb{o}  \\
		\vb{o} & \vb{o} & \curl \bm{\theta}_1 &  & \vb{o} & \vb{o} & \curl\bm{\theta}_6  
	\end{bmatrix} \in \mathbb{R}^{9 \times 18} \, .
\end{align}

Consequently, the element stiffness matrix reads
\begin{align}
	\bm{K}_\text{n\'ed} = \int_\Omega (\bm{\Theta}^T\mathbb{S}\, \bm{\Theta} + \Curl (\bm{\Theta})^T \, \mathbb{S} \, \Curl \bm{\Theta}) \, \dd X \in \mathbb{R}^{18 \times 18} \, .
\end{align}
The N\'ed\'elec element belongs to the subspace 
$\Ned^0_I \subset \Hc{}$. The corresponding interpolant yields the following error estimate
\begin{align}
    \| \vb{v} - \Pi_{\Ned} \vb{v} \|_{\Le} \preceq h \,  |\vb{v}|_{\mathit{H}^1} \, ,
\end{align}
where $\vb{v}$ is the exact solution, $\Pi_{\Ned}$ is the N\'ed\'elec interpolant and $h$ is a measure of the element's size.

\subsection{The $\Hsc{}$-element}
Due to the following inclusion the two former finite element formulations are possible candidates for computations in $\Hsc{}$. However, as denoted by the trace in \cref{eq:symcurltr} and shown in \cite{Lewintan2021}, the $\Hsc{}$ space is strictly larger
\begin{align}
	[\Hone]^{3 \times 3} \subset [\Hc{}]^3 \subset \Hsc{} \, .
\end{align}
Consequently, there can exist solutions that belong to $\Hsc{}$ while not belonging to $[\Hc{}]^3$, as demonstrated in \cref{sec:dis}.
For the design of an $\Hsc{}$-conforming finite element \cite{sander2021conforming}, the trace must vanish on the interface of neighbouring elements. In other words, a function $\Pm$ belongs to $\Hsc{}$ iff
\begin{align}
	\jump{\tr_{\Hsc{}}(\Pm)} \at_{\Xi_i} = 0 \qquad \forall \, \Xi_i  \text{ in } \Omega \, , 
\end{align}
where $\jump{\cdot}$ represents the jump operator and $\Xi_i$ is an interface between neighbouring elements.
The trace condition can be reconstructed as 
\begin{align}
	&\sym[\Pm \Anti(\bm{\nu})^T] = 0 \quad \Longleftrightarrow \quad \langle \Pm \Anti(\bm{\nu})^T , \, \bm{S}_i \rangle = 0  \qquad \forall \, \bm{S}_i \in \Sym(3) \, . 
	\label{eq:newtrace}
\end{align}
A simple basis for $\Sym(3)$ is given by the symmetric decomposition of a matrix filled with ones (${\bm{A} = \sum_i \sum_j 1\,\vb{e}_i \otimes \vb{e}_j}$) 
\begin{align}
	\Sym(3) = \spa \{ \vb{e}_1 \otimes \vb{e}_1 , \; \vb{e}_2 \otimes \vb{e}_2 , \; \vb{e}_3 \otimes \vb{e}_3 , \; \sym(\vb{e}_1 \otimes \vb{e}_2) , \; \sym(\vb{e}_1 \otimes \vb{e}_3) , \; \sym(\vb{e}_2 \otimes \vb{e}_3) \} \, .
\end{align}
Accordingly, any affine transformation $\bm{F}$ that maps to linearly independent vectors allows to generalize the basis
\begin{align}
	&\Sym(3) = \spa \{ \bm{F} \, \vb{e}_1 \otimes \vb{e}_1 \bm{F}^T  , \dots , \; \sym(\bm{F}\, \vb{e}_2 \otimes \vb{e}_3 \bm{F}^T) \} \, , && \bm{F} \in \GL(3) \, .
\end{align} 
Clearly, the transformation leaves the symmetry invariant. 
Defining $\vb{a} \in \mathbb{R}^3$, $\vb{b} \in \mathbb{R}^3$ such that \\ ${\bm{F} = \begin{bmatrix}
	\vb{a} & \vb{b} & \bm{\nu}
\end{bmatrix} \in \GL(3)}$, where $\bm{\nu}$ is a normal on an element's face, we observe that five conditions instead of six suffice to assert \cref{eq:newtrace} since
\begin{align}
	\langle \Pm \Anti(\bm{\nu})^T , \, \bm{\nu} \otimes \bm{\nu} \rangle = \langle \Pm , \,  \bm{\nu} \otimes \bm{\nu} \Anti(\bm{\nu}) \rangle = \langle \Pm , \,  \bm{\nu} \otimes \bm{\nu} \times \bm{\nu} \rangle = 0 \, ,
\end{align}
is always satisfied. In other words, five conditions control whether $\Pm$ is point-wise $\Hsc{}$-conforming on a face
\begin{align}
	&\langle \Pm \Anti(\bm{\nu})^T , \, \vb{a} \otimes \vb{a} \rangle = 0 \, , && \langle \Pm \Anti(\bm{\nu})^T , \, \vb{b} \otimes \vb{b} \rangle = 0 \, , && \langle \Pm \Anti(\bm{\nu})^T , \, \sym(\vb{a} \otimes \vb{b}) \rangle = 0 \, , \notag \\
	&\langle \Pm \Anti(\bm{\nu})^T , \, \sym(\vb{a} \otimes \bm{\nu}) \rangle = 0 \, , && \langle \Pm \Anti(\bm{\nu})^T , \, \sym(\vb{b} \otimes \bm{\nu}) \rangle = 0 \, . 
	\label{eq:cond}
\end{align}
Yet, these are not enough conditions to fully identify an element of $\mathbb{R}^{3 \times 3}$. The remaining conditions can be defined as 
\begin{align}
	&\langle \Pm , \,  \vb{b} \otimes (\vb{a} \times \bm{\nu})  - \vb{a} \otimes (\vb{b} \times \bm{\nu}) \rangle = 0 \, , && 
	\langle \Pm , \, \vb{a} \otimes \bm{\nu} \rangle = 0 \, , \notag \\
	&\langle \Pm , \, \vb{b} \otimes\bm{\nu} \rangle = 0 \, , && 
	\langle \Pm , \, \bm{\nu} \otimes \bm{\nu} \rangle = 0  \, , 
\end{align}
which, together with the trace conditions, determine the element uniquely.
Note that the latter identities are all in $\ker(\tr_{\Hsc{}})$ due to $\bm{\nu} \times \bm{\nu} = 0$ and are therefore linearly independent of the trace conditions. To see their linear independence from one another, simply map back to the Cartesian basis with $\bm{F}^{-1}$.

The conditions in \cref{eq:cond} ensure conformity on a face. On an edge two faces meet and the conditions must be satisfied for both. On each edge we define $\vb{a} = \bm{\tau}$, $\vb{b} = \bm{\gamma}_i = \bm{\tau} \times \bm{\nu}_i$, where $\bm{\tau}$ is the edge tangent, $\bm{\nu}_i$ is the normal on face $i$ and $\bm{\gamma}_i$ is the corresponding conormal.
Considering the following identity
\begin{align}
	\bm{\nu}_1 \times \bm{\gamma}_1 = \bm{\nu}_2 \times \bm{\gamma}_2 = \bm{\tau} \, ,
\end{align} 
the conditions for conformity can be reformulated as
\begin{align}
	&\langle \Pm , \, \bm{\tau} \otimes \bm{\gamma}_i \rangle = 0 \, , && \langle \Pm , \, \bm{\tau} \otimes \bm{\tau} 
	+ \bm{\gamma}_i \otimes \bm{\gamma}_i \rangle = 0 \, ,
	&& \langle \Pm , \, \bm{\nu}_i \otimes \bm{\gamma}_i \rangle = 0 \, , && i = \{1, \, 2\} \,  , \notag \\
	 &\langle \Pm , \, \bm{\tau} \otimes \bm{\rho}_1 \rangle = 0 \, , && \langle \Pm , \, \bm{\tau} \otimes \bm{\rho}_2  \rangle = 0 \, ,
\end{align}
where $\bm{\rho}_1$ and $\bm{\rho}_2$ are two arbitrary vectors spanning the surface orthogonal to $\bm{\tau}$.
The upper conditions allow to determine eight terms of $\Pm$. One term remains to be determined using 
\begin{align}
	\langle \Pm , \, \bm{\tau} \otimes \bm{\tau} \rangle = 0 \, ,
\end{align} 
which controls multiples of the identity matrix. The conditions are linearly independent as long as $\bm{\nu}_1 \nparallel \bm{\nu}_2$.

The same methodology may be used to construct vertex conformity and uniqueness conditions. However, for an unstructured mesh it is unclear how to do this in a way that is independent of the geometry. Consequently, at the vertices we impose full continuity with the exception of the identity 
\begin{align}
	&P_{12}, \, P_{13}, \, P_{21}, \, P_{23}, \, P_{31}, \, P_{32} = 0 \, , && P_{11} - P_{22} = 0 \, , && P_{22} - P_{33} = 0 \, , \notag \\ &\tr \Pm = \langle \Pm , \, \one \rangle = P_{11} + P_{22} + P_{33} = 0 \, ,
	\label{eq:dofs}
\end{align}  
where the conditions in the upper row assert conformity and the trace controls the identity matrix.

Unlike for the previous elements that were constructed in the reference domain, the construction of the $\Hsc{}$-element is done directly on the grid using the degrees of freedom from \cref{eq:dofs}. 
We reformulate the degrees of freedom for the linear element as vector operators
\begin{align}
	&P_{11} - P_{22} : \quad \vb{l}_1 = \begin{bmatrix}
		1 & 0 & 0 & 0 & -1 & 0 & 0 & 0 & 0
	\end{bmatrix} \, , && P_{22} - P_{33} : \quad 
    \vb{l}_5 = \begin{bmatrix}
    	0 & 0 & 0 & 0 & 1 & 0 & 0 & 0 & -1
    \end{bmatrix} \, , \notag \\ 
	&\tr \Pm: \quad \vb{l}_9 = \begin{bmatrix}
	1 & 0 & 0 & 0 & 1 & 0 & 0 & 0 & 1
    \end{bmatrix} \, , && P_{ij} \quad \text{s.t.} \quad i\neq j: \quad
     \vb{l}_i = \vb{e}_i^T \, , \quad i \in \{2,3,4,6,7,8\} \, ,
\end{align}
where $\vb{e}_i$ is the unit vector in $\mathbb{R}^9$.
Consequently, the collection of the degrees of freedom can be defined as the operator matrix 
\begin{align}
	\bm{L} = \begin{bmatrix}
		\vb{l}_1^T & \vb{l}_2^T & \vb{l}_3^T & \vb{l}_4^T & \vb{l}_5^T & \vb{l}_6^T & \vb{l}_7^T & \vb{l}_8^T & \vb{l}_9^T
	\end{bmatrix}^T \in \mathbb{R}^{9 \times 9 } \, .
\end{align}
The ansatz matrix is given by the monomial basis
\begin{align}
	\bm{X} = \begin{bmatrix}
		1 \one_9 & x \one_9 &  y \one_9 &  z \one_9 
	\end{bmatrix} \in \mathbb{R}^{9 \times 36} \, , && \one_9 \in \mathbb{R}^{9 \times 9} \, .
\end{align}
In order to construct the corresponding base functions, the product $\bm{L} \, \bm{X}$ is evaluated at each node of a finite element
\begin{align}
	\bm{C}^{-1} = \begin{bmatrix}
		(\bm{L} \, \bm{X})^T \at_{n_1} & (\bm{L} \, \bm{X})^T \at_{n_2} & (\bm{L} \, \bm{X})^T \at_{n_3} & (\bm{L} \, \bm{X})^T \at_{n_4}
	\end{bmatrix}^T \in \mathbb{R}^{36 \times 36} \, ,
\end{align}
where every column of $\bm{C}$ gives the constant factors for one base function. Therefore, the resulting local basis is given by 
\begin{align}
	&\bm{B} = \bm{X} \, \bm{C} \in \mathbb{R}^{9 \times 36} \, , && \bm{P} = \bigcup_{e=1}^{n_{\text{elements}}} \bm{B} \,  \bm{P}_e \, .
\end{align}
In order to find the Curl of $\Pm$, we redefine the operator in Voigt notation
\begin{align}
		\Curl \bm{P} = \begin{bmatrix}
			\Anti (\nabla) & \bm{O} & \bm{O} \\
			\bm{O} &\Anti (\nabla) & \bm{O} \\
			\bm{O} & \bm{O}& \Anti (\nabla) 
		\end{bmatrix} \, \bm{P} \, , && 
		\Anti(\nabla) = \begin{bmatrix}
			0 & -\partial_3 & \partial_2 \\
			\partial_3 & 0  & -\partial_1 \\
			-\partial_2 & \partial_1  & 0 
		\end{bmatrix} \, ,
		&&\bm{O} = \begin{bmatrix}
			0 & 0 & 0 \\
			0 & 0 & 0 \\
			0 & 0 & 0 
		\end{bmatrix} \, .
\end{align}
Consequently, the Curl of the base functions is given by
\begin{align}
	\Curl \bm{B} = (\Curl \bm{X}) \, \bm{C} \, .
\end{align}
The element stiffness matrix can now be defined as
\begin{align}
	\bm{K}_\text{sym} = \int_\Omega (\bm{B}^T\mathbb{S}\, \bm{B} + \Curl (\bm{B})^T \, \mathbb{S} \, \Curl \bm{B}) \, \dd X \in \mathbb{R}^{36 \times 36} \, .
\end{align}

\begin{remark}
	At the construction of the local-global map, eight of the nine degrees of freedom on a node are shared by neighbouring elements. The ninth degree is built locally for each element and is not connected over the element's boundaries. With this characteristic, the lower continuity of the space is made possible.
\end{remark}
The element belongs to $\San^1 \subset \Hsc{}{}$ and its interpolant produces the standard Lagrangian error estimate (see \cite{sander2021conforming}) 
\begin{align}
    \| \bm{T} - \Pi_{\San} \bm{T} \|_{\Le} \preceq h^2 \,  |\bm{T}|_{\mathit{H}^2} \, ,
\end{align}
where $\bm{T}$ is the exact solution tensor, $\Pi_{\San}$ is the interpolant and $h$ is a measure of the element's size.

\section{Numerical examples}\label{seq:eg}
\begin{figure}
	\centering
	\begin{subfigure}{0.3\linewidth}
		\includegraphics[width=1.0\linewidth]{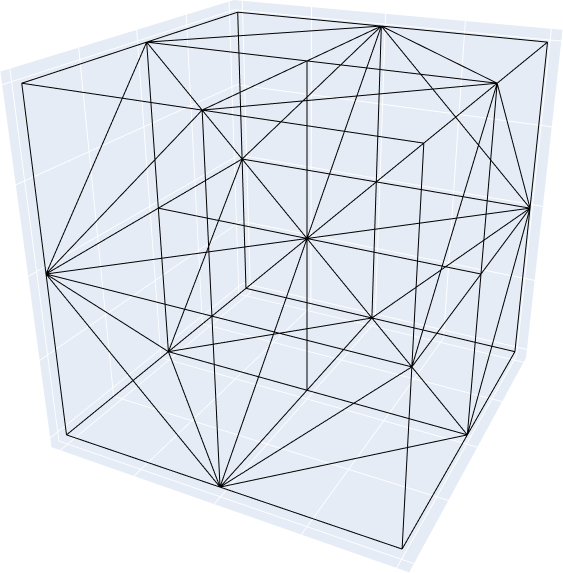}
		\caption{40 elements} 
	\end{subfigure}
    \begin{subfigure}{0.3\linewidth}
    	\includegraphics[width=1.0\linewidth]{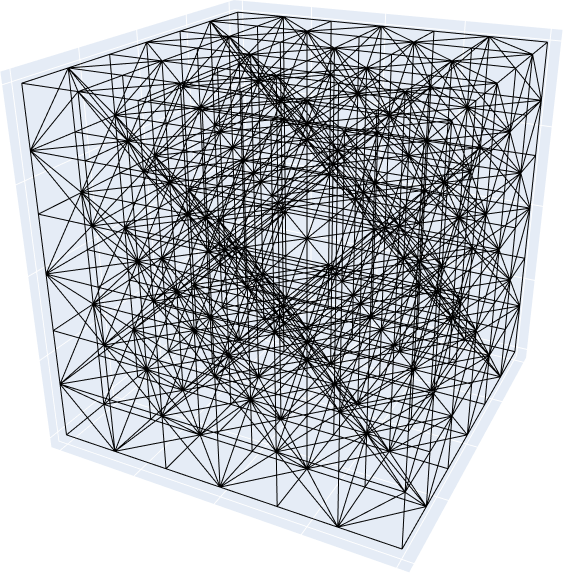}
    	\caption{1080 elements}
    \end{subfigure}
    \begin{subfigure}{0.3\linewidth}
    	\includegraphics[width=1.0\linewidth]{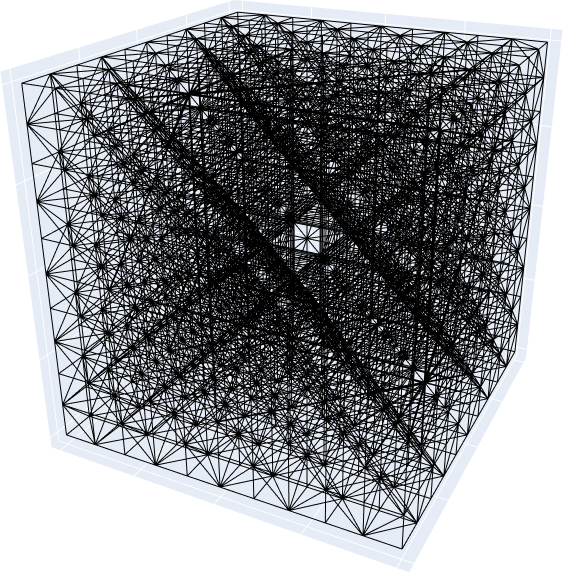}
    	\caption{5000 elements}
    \end{subfigure}
    \caption{Finite element meshes from $40$ to $5000$ elements for the domain $\Omega = [-1, \, 1]^3$.}
    \label{fig:mesh}
\end{figure}

In the following we test the convergence of the three finite element formulations against various simple analytical solutions. 
For all formulations we make use of the cube $\Omega = [-1,\, 1]^3$. The finite element meshes range from $40$ to $5000$ elements, see \cref{fig:mesh}.

We measure convergence in the Lebesgue $\Le$-norm
\begin{align}
	\|\widetilde{\Pm} - \Pm \|_{\Le} = \sqrt{\int_\Omega \| \widetilde{\Pm} - \Pm \|_F^2 \dd X } \, ,
\end{align}
and the $\Hsc{}$-norm
\begin{align}
    \|\widetilde{\Pm} - \Pm \|_{\Hsc{}} = \sqrt{\|\widetilde{\Pm} - \Pm \|^2_{\Le} + \int_\Omega \| \sym \Curl (\widetilde{\Pm} - \Pm) \|_F^2 \dd X}
    \, .
\end{align}

\subsection{Smooth solution field}
\begin{figure}
	\centering
	\begin{subfigure}{0.48\linewidth}
		\begin{tikzpicture}
			\begin{loglogaxis}[
				/pgf/number format/1000 sep={},
				axis lines = left,
				xlabel={Degrees of freedom},
				ylabel={$\| \widetilde{\bm{P}} - \bm{P}\|_{\Le}$},
				xmin=100, xmax=1e5,
				ymin=0.01, ymax=100,
				xtick={100,1000,10000,1e5},
				ytick={0.01, 0.1, 1, 10, 100},
				legend pos = north east,
				ymajorgrids=true,
				grid style=dotted,
				]
				\addplot[
				color=black,
				mark=triangle,
				]
				coordinates {
					( 243 , 2.001268188377093 )
					( 1125 , 0.49012270906357147 )
					( 3087 , 0.21386047076536738 )
					( 6561 , 0.11865193672577146 )
					( 11979 , 0.07523348381595904 )
				};
				\addlegendentry{$[\Hone{}]^{3\times3}$}
				\addplot[
				color=teal,
				mark=o,
				]
				coordinates {
					( 270 , 2.204563785578984 )
					( 1620 , 1.049169285065994 )
					( 4914 , 0.6867814954169372 )
					( 11016 , 0.5106343028127667 )
					( 20790 , 0.406501032452444 )
				};
				\addlegendentry{$[\Hc{}]^3$}
				\addplot[
				color=blue,
				mark=square,
				]
				coordinates {
					( 376 , 1.785569726956325 )
					( 2280 , 0.4532601127990592 )
					( 7064 , 0.2007905656715176 )
					( 16072 , 0.11246696209877236 )
					( 30648 , 0.07178333006846843 )
				};
				\addlegendentry{$\Hsc{}$}
				\addplot[dashed,color=black, mark=none]
				coordinates {
					(500, 4)
					(25000, 1.0857670466379628)
				};
				\addplot[dashed,color=black, mark=none]
				coordinates {
					(200, 0.2)
					(10000, 0.014736125994561549)
				};
			\end{loglogaxis}
			\draw (3.,3.9) node[anchor=north west]{$\mathcal{O}(h)$};
			\draw (2.5,1.6) node[anchor=north west]{$\mathcal{O}(h^2)$};
		\end{tikzpicture}
	\end{subfigure}
	\begin{subfigure}{0.48\linewidth}
		\begin{tikzpicture}
			\begin{loglogaxis}[
				/pgf/number format/1000 sep={},
				axis lines = left,
				xlabel={Degrees of freedom},
				ylabel={$\|\widetilde{\bm{P}} - \bm{P}\|_{\Hsc{}}$},
				xmin=100, xmax=1e5,
				ymin=0.01, ymax=100,
				xtick={100,1000,10000,1e5},
				ytick={0.01, 0.1, 1, 10, 100},
				legend pos = north east,
				ymajorgrids=true,
				grid style=dotted,
				]
				\addplot[
				color=black,
				mark=triangle,
				]
				coordinates {
					( 243 , 5.255325784707402 )
					( 1125 , 2.647609048932859 )
					( 3087 , 1.7561053613632522 )
					( 6561 , 1.3114997610793833 )
					( 11979 , 1.045894715062534 )
				};
				\addlegendentry{$[\Hone{}]^{3\times3}$}
				\addplot[
				color=teal,
				mark=o,
				]
				coordinates {
					( 270 , 5.048746406564437 )
					( 1620 , 2.619432424035243 )
					( 4914 , 1.7459776913749037 )
					( 11016 , 1.3071427138000666 )
					( 20790 , 1.0440954201601618 )
				};
				\addlegendentry{$[\Hc{}]^3$}
				\addplot[
				color=blue,
				mark=square,
				]
				coordinates {
					( 376 , 5.177029456302448 )
					( 2280 , 2.641033535551179 )
					( 7064 , 1.7545616519562532 )
					( 16072 , 1.310954674573544 )
					( 30648 , 1.0456522004875117 )
				};
				\addlegendentry{$\Hsc{}$}
				\addplot[dashed,color=black, mark=none]
				coordinates {
					(500, 0.5)
					(25000, 0.13572088082974534)
				};
			\end{loglogaxis}
			\draw (3.,2.6) node[anchor=north west]{$\mathcal{O}(h)$};
		\end{tikzpicture}
	\end{subfigure}
	\caption{Convergence rates for the rotational benchmark.}
	\label{fig:vanish}
\end{figure}
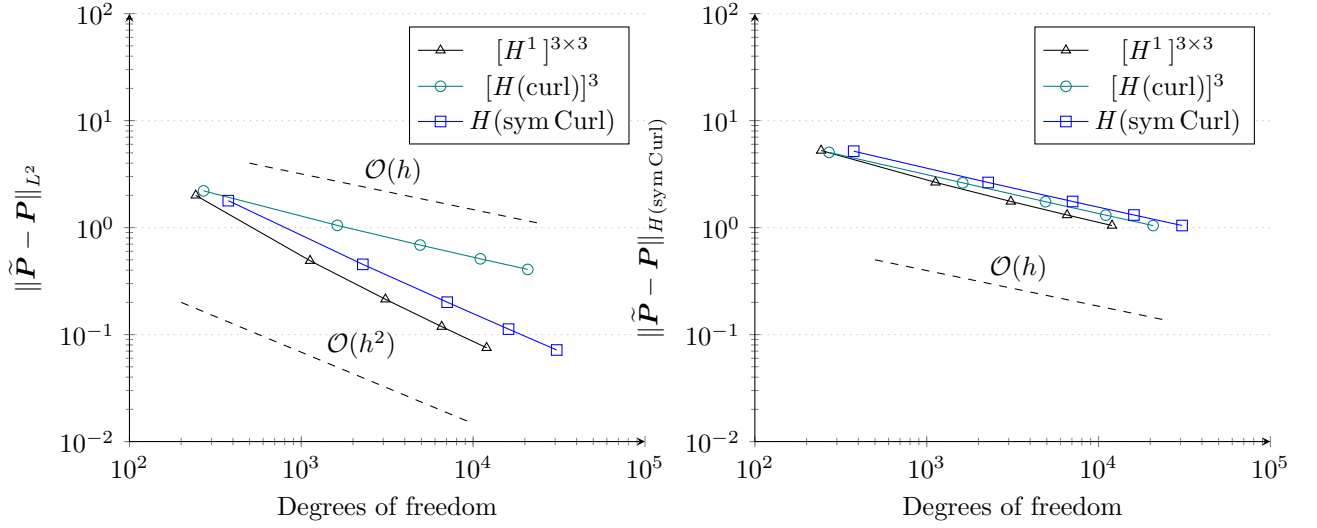

In the first example we set the microdistortion to
\begin{figure}
	\centering
    \begin{subfigure}{0.3\linewidth}
    	\includegraphics[width=1.0\linewidth]{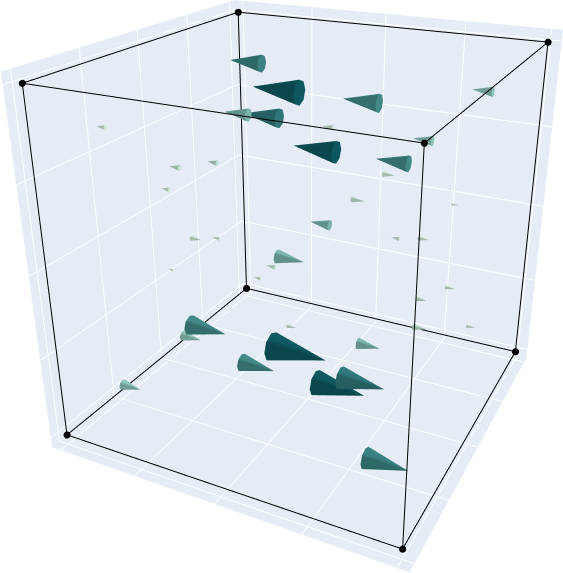}
    	\caption{40 elements} 
    \end{subfigure}
    \begin{subfigure}{0.3\linewidth}
    	\includegraphics[width=1.0\linewidth]{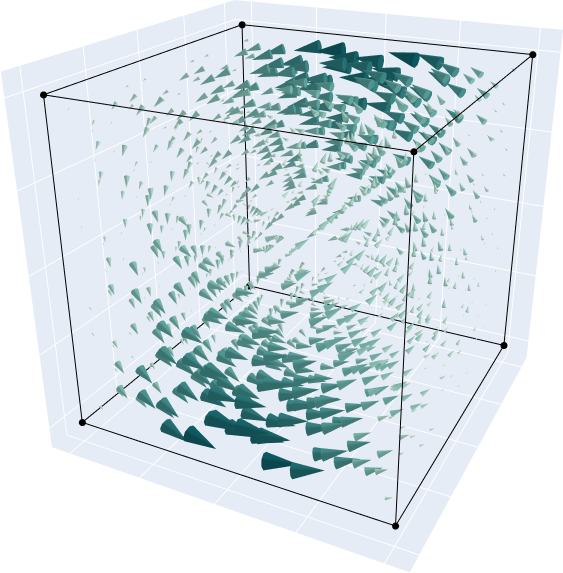}
    	\caption{1080 elements}
    \end{subfigure}
    \begin{subfigure}{0.3\linewidth}
    	\includegraphics[width=1.0\linewidth]{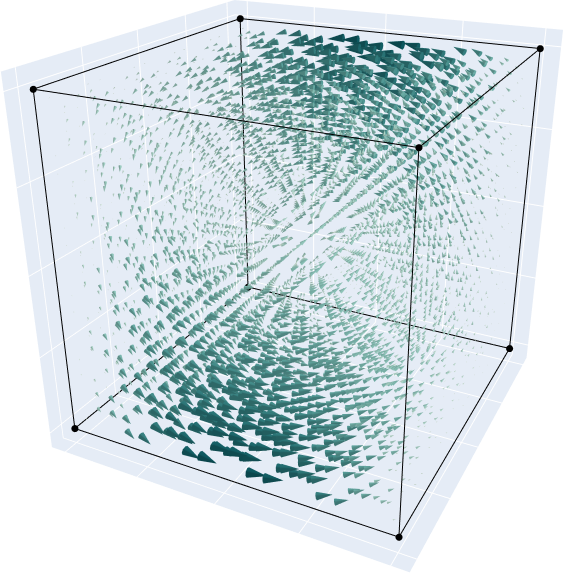}
    	\caption{5000 elements}
    \end{subfigure}
	\caption{Vortex field for various discretizations of the $\Hsc{}$-element displaying the first row of $\Pm$.}
	\label{fig:vanishpics}
\end{figure}
\begin{align}
	\widetilde{\Pm} = (1-x)(1+x) \begin{bmatrix}
		-y -z & x & x \\
		-y -z & x & x \\
		-y -z & x & x 
	\end{bmatrix} \, .
\end{align} 
Effectively, the microdistortion is now a vortex field that must vanish at $x = -1$ and $x = 1$, see \cref{fig:vanishpics}. 
Since the microdistortion is clearly continuous, it is an element of $[\Hone]^{3 \times 3}$.
The microdistortion field gives rise to the micro-moment
\begin{align}
	\bm{M} =  \begin{matrix}
		\begin{bmatrix}x^{2} y +  x^{2} z -  y -  z & 0.5(x^{2} y-  x^{3} + x^{2} z + 9 x - y - z)  & 0.5(x^{2} y -  x^{3}   + x^{2} z + 9 x -  y -  z) \\0.5(x^{2} y -  x^{3} + x^{2} z +
			x - y - z) & - x^{3} + x & - x^{3} + 9 x\\0.5(x^{2} y-  x^{3}  + x^{2} z + x - y - z) & - x^{3} +  9x & -  x^
			{3} + x\end{bmatrix}
		& \, .
	\end{matrix} 
\end{align}
The convergence rates are given in \cref{fig:vanish}.
In the $\Le$-norm, the N\'ed\'elec element converges linearly, whereas the Lagrangian element and the $\Hsc{}$-element converge quadratically, as expected for smooth fields.
All three formulations converge linearly in the $\Hsc{}$-norm.

\subsection{Discontinuous normal trace}\label{sec:trcurl}
\begin{figure}
	\centering
	\begin{subfigure}{0.48\linewidth}
		\begin{tikzpicture}
			\begin{loglogaxis}[
				/pgf/number format/1000 sep={},
				axis lines = left,
				xlabel={Degrees of freedom },
				ylabel={$\| \widetilde{\bm{P}} - \bm{P}\|_{\Le}$},
				xmin=100, xmax=1e5,
				ymin=0.01, ymax=100,
				xtick={100,1000,10000,1e5},
				ytick={0.01, 0.1, 1, 10, 100},
				legend pos = north east,
				ymajorgrids=true,
				grid style=dotted,
				]
				\addplot[
				color=black,
				mark=triangle,
				]
				coordinates {
					( 243 , 1.04741581514085 )
					( 1125 , 0.7297484247203386 )
					( 3087 , 0.5770578809456791 )
					( 6561 , 0.4956099951705984 )
					( 11979 , 0.43838561334086934 )
				};
				\addlegendentry{$[\Hone{}]^{3\times3}$}
				\addplot[
				color=teal,
				mark=o,
				]
				coordinates {
					( 270 , 4.643863056372025e-16 )
					( 1620 , 7.240670470274879e-15 )
					( 4914 , 2.0565084669697132e-14 )
					( 11016 , 4.456267283748485e-14 )
					( 20790 , 9.319243482207808e-14 )
				};
				\addlegendentry{$[\Hc{}]^3$}
				\addplot[
				color=blue,
				mark=square,
				]
				coordinates {
					( 376 , 0.8915479053785873 )
					( 2280 , 0.6281726161959351 )
					( 7064 , 0.5043109872632786 )
					( 16072 , 0.4339772126762664 )
					( 30648 , 0.3856033692913387 )
				};
				\addlegendentry{$\Hsc{}$}
				\addplot[dashed,color=black, mark=none]
				coordinates {
					(500, 0.2)
					(25000, 0.10420014619173827)
				};
			\end{loglogaxis}
			\draw (3.,2.25) node[anchor=north west]{$\mathcal{O}(\sqrt{h})$};
		\end{tikzpicture}
	\end{subfigure}
	\begin{subfigure}{0.48\linewidth}
		\begin{tikzpicture}
			\begin{loglogaxis}[
				/pgf/number format/1000 sep={},
				axis lines = left,
				xlabel={Degrees of freedom },
				ylabel={$\|\widetilde{\bm{P}} - \bm{P}\|_{\Hsc{}}$ },
				xmin=100, xmax=1e5,
				ymin=0.01, ymax=100,
				xtick={100,1000,10000,1e5},
				ytick={0.01, 0.1, 1, 10, 100},
				legend pos = north east,
				ymajorgrids=true,
				grid style=dotted,
				]
				\addplot[
				color=black,
				mark=triangle,
				]
				coordinates {
					( 243 , 1.0674239735896642 )
					( 1125 , 0.7449439515189085 )
					( 3087 , 0.592632555745497 )
					( 6561 , 0.5096126004950042 )
					( 11979 , 0.45140828429701163 )
				};
				\addlegendentry{$[\Hone{}]^{3\times3}$}
				\addplot[
				color=teal,
				mark=o,
				]
				coordinates {
					( 270 , 7.448764101789542e-16 )
					( 1620 , 7.91846627548546e-15 )
					( 4914 , 2.146573552245785e-14 )
					( 11016 , 4.652851648848629e-14 )
					( 20790 , 9.527585396753145e-14 )
				};
				\addlegendentry{$[\Hc{}]^3$}
				\addplot[
				color=blue,
				mark=square,
				]
				coordinates {
					( 376 , 0.9149708832371244 )
					( 2280 , 0.6457627763050725 )
					( 7064 , 0.5220604563097896 )
					( 16072 , 0.4499021631338771 )
					( 30648 , 0.4003466642255849 )
				};
				\addlegendentry{$\Hsc{}$}
				\addplot[dashed,color=black, mark=none]
				coordinates {
					(500, 0.2)
					(25000, 0.10420014619173827)
				};
			\end{loglogaxis}
			\draw (3.,2.25) node[anchor=north west]{$\mathcal{O}(\sqrt{h})$};
		\end{tikzpicture}
	\end{subfigure}
	\caption{Convergence for the discontinuous normal projection benchmark.}
	\label{fig:trcurl}
\end{figure}
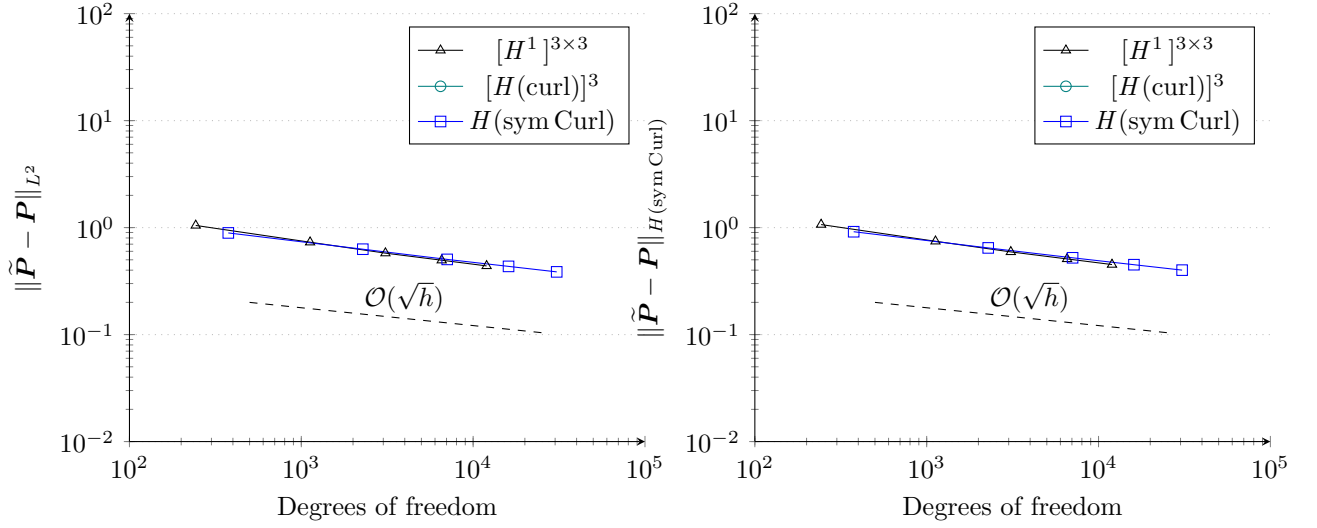

\begin{figure}
	\centering
	\begin{subfigure}{0.3\linewidth}
		\includegraphics[width=1.0\linewidth]{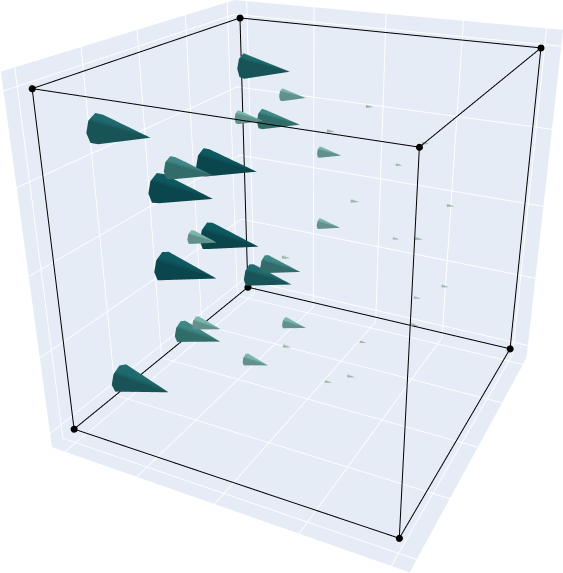}
		\caption{40 elements} 
	\end{subfigure}
	\begin{subfigure}{0.3\linewidth}
		\includegraphics[width=1.0\linewidth]{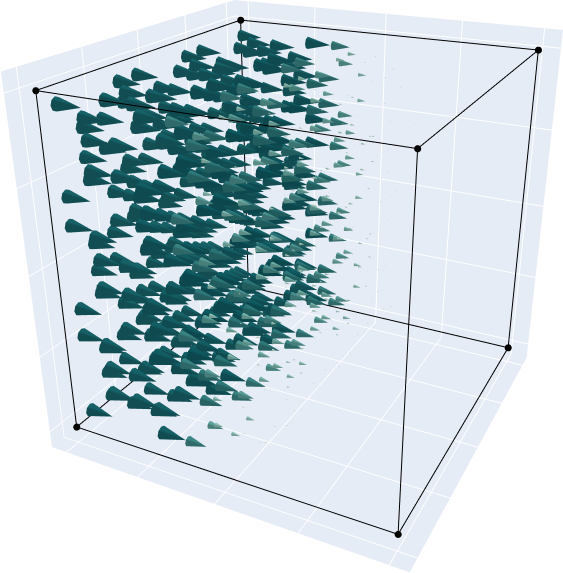}
		\caption{1080 elements}
	\end{subfigure}
	\begin{subfigure}{0.3\linewidth}
		\includegraphics[width=1.0\linewidth]{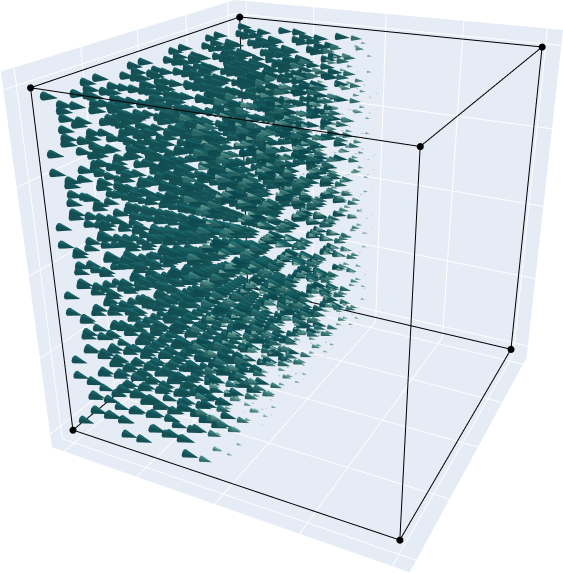}
		\caption{5000 elements}
	\end{subfigure}
	\caption{Discontinuous normal projection for various discretizations of the $\Hone$-element displaying the first row of $\Pm$.}
	\label{fig:normallag}
\end{figure}

\begin{figure}
	\centering
	\begin{subfigure}{0.3\linewidth}
		\includegraphics[width=1.0\linewidth]{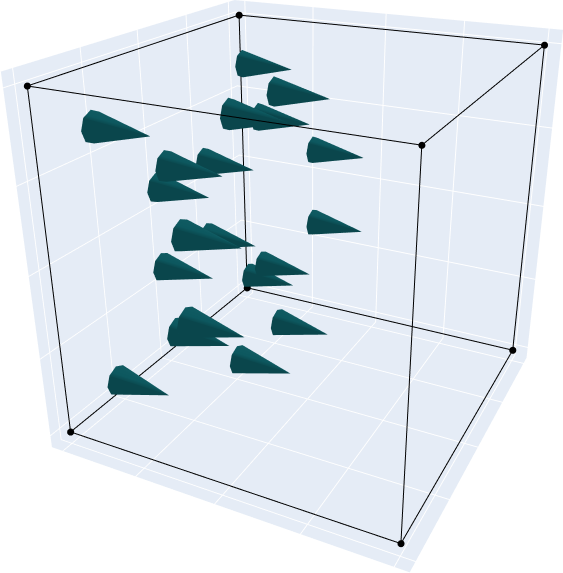}
		\caption{40 elements} 
	\end{subfigure}
	\begin{subfigure}{0.3\linewidth}
		\includegraphics[width=1.0\linewidth]{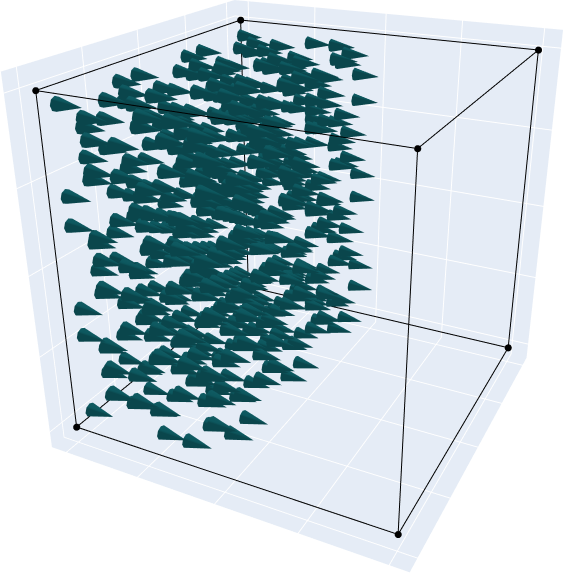}
		\caption{1080 elements}
	\end{subfigure}
	\begin{subfigure}{0.3\linewidth}
		\includegraphics[width=1.0\linewidth]{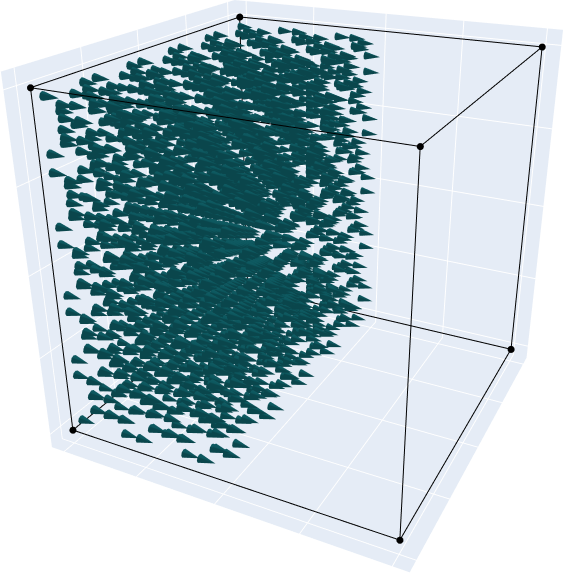}
		\caption{5000 elements}
	\end{subfigure}
	\caption{Discontinuous normal projection for various discretizations of the $\Hc{}$-element displaying the first row of $\Pm$.}
	\label{fig:normalned}
\end{figure}

\begin{figure}
	\centering
	\begin{subfigure}{0.3\linewidth}
		\includegraphics[width=1.0\linewidth]{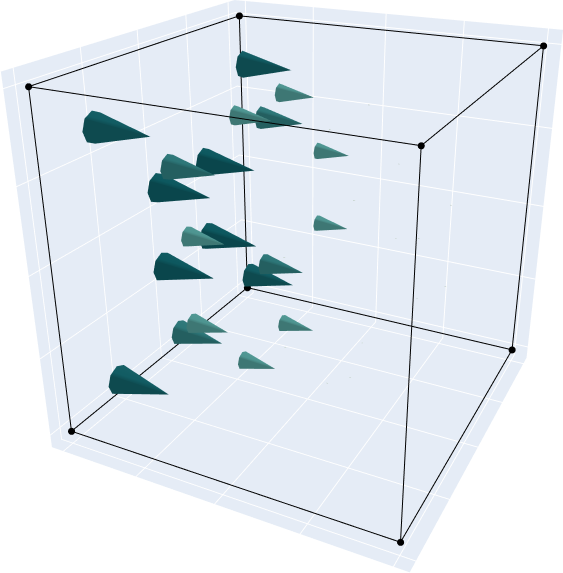}
		\caption{40 elements} 
	\end{subfigure}
	\begin{subfigure}{0.3\linewidth}
		\includegraphics[width=1.0\linewidth]{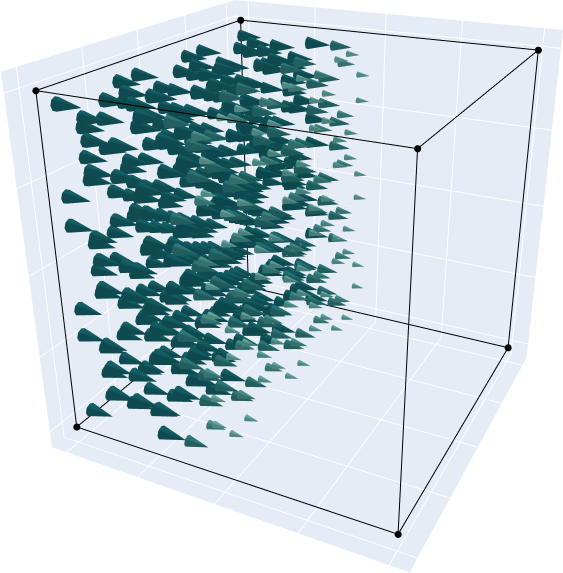}
		\caption{1080 elements}
	\end{subfigure}
	\begin{subfigure}{0.3\linewidth}
		\includegraphics[width=1.0\linewidth]{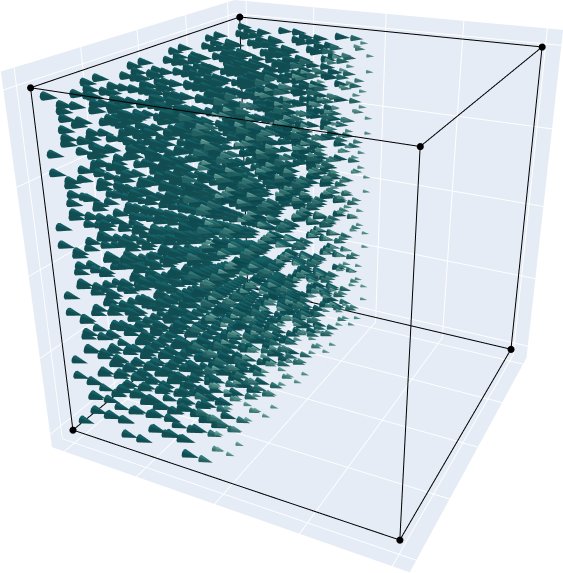}
		\caption{5000 elements}
	\end{subfigure}
	\caption{Discontinuous normal projection for various discretizations of the $\Hsc{}$-element displaying the first row of $\Pm$.}
	\label{fig:normalsc}
\end{figure}

In the following benchmark we define the discontinuous microdistortion field
\begin{align}
	\widetilde{\Pm} = \left \{  \begin{matrix}
		\begin{bmatrix}
		1 & 0 & 0 \\
		0 & 0 & 0 \\
		0 & 0 & 0
		\end{bmatrix} & \quad \text{for} & x < 0  \\[4ex]
		0 & \quad \text{otherwise} &
	\end{matrix}  \right. \, .
\end{align}
The latter field belongs to $[\Hc{}]^3$ but not to $[\Hone]^{3 \times 3}$. This is because only the normal projection of the microdistortion $\Pm$ is discontinuous
with regard to the unit vector $\vb{e}_1 = \begin{bmatrix}
1 & 0 & 0
\end{bmatrix}^T$
\begin{align}
    \tr_{\Hc{}} \Pm \at_{\Gamma_1} = \Pm \Anti(\vb{e}_1)^T = 0 \, , && \Gamma_1 \perp \vb{e}_1 \, .
\end{align}
The corresponding micro-moments are clearly
\begin{align}
	\bm{M} = \widetilde{\Pm} \, .
\end{align}
The approximation captured by the various discretizations is depicted in \cref{fig:normallag}, \cref{fig:normalned}
and \cref{fig:normalsc}. The noise (small vectors in $x > 0$) in the solution is apparent for the Lagrangian- and $\Hsc{}$-formulations. 
The convergence rates are given in \cref{fig:trcurl}. The N\'ed\'elec element finds the analytical solution immediately for all discretizations. Both the Lagrangian formulation and the $\Hsc{}$-formulation converge sub-optimally. In case of the $\Hsc{}$-formulation this is due to the higher continuity imposed at the vertices.
Despite appearing similar, the values of the convergence in the $\Le$-norm are not equal to the values of the convergence in the $\Hsc{}$-norm. However, the difference is small.

\subsection{Discontinuous identity}\label{sec:dis}
In the last benchmark we test the behaviour of a discontinuous identity microdistortion field
\begin{align}
	\widetilde{\Pm} = \left \{  \begin{matrix}
		\one & \quad \text{for} & x < 0  \\
		0 & \quad \text{otherwise} &
	\end{matrix}  \right. \, ,
\end{align}
for which the micro-moment reads
\begin{align}
	\bm{M} = \widetilde{\Pm} \, .
\end{align}
Clearly the solution is in $\Hsc{}$ but not in $[\Hc{}]^3$ or $[\Hone]^{3 \times 3}$ since
\begin{align}
    \tr_{\Hc{}} \Pm \at_{\Gamma_1} &= \left \{  \begin{matrix}
		\Anti(\vb{e}_1)^T & \quad \text{for} & x < 0  \\
		0 & \quad \text{otherwise} &
	\end{matrix}  \right.  \, ,  \notag \\
	\tr_{\Hsc{}} \Pm \at_{\Gamma_1} &= 0 \, ,
\end{align}
where $\vb{e}_1 = \begin{bmatrix}
1 & 0 & 0
\end{bmatrix}^T$ and $\Gamma_1 \perp \vb{e}_1$.

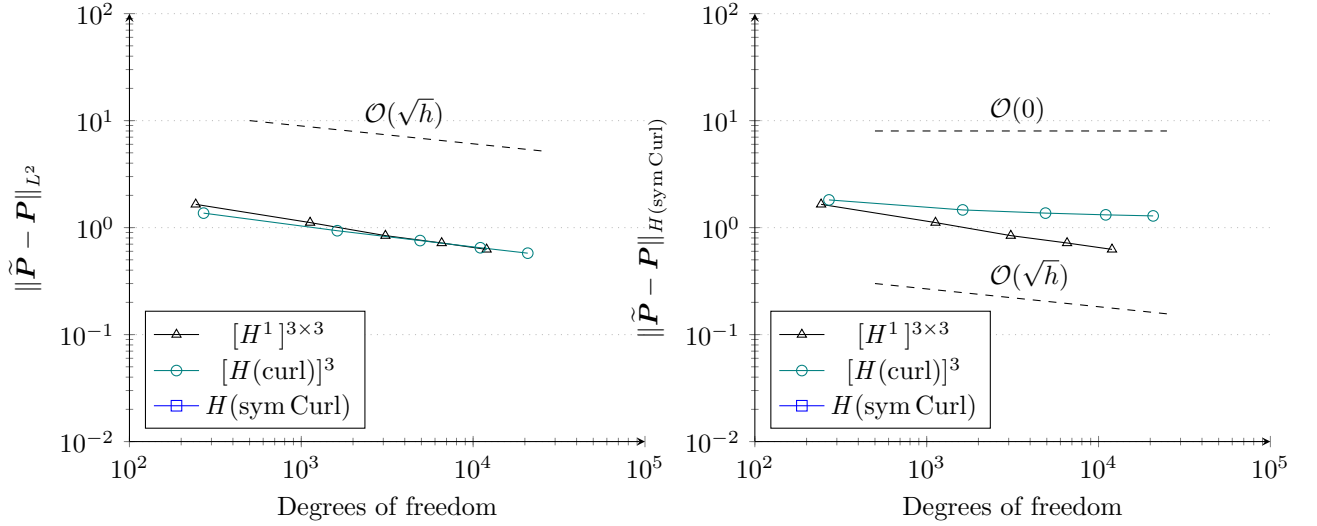
\begin{figure}
	\centering
	\begin{subfigure}{0.48\linewidth}
		\begin{tikzpicture}
			\begin{loglogaxis}[
				/pgf/number format/1000 sep={},
				axis lines = left,
				xlabel={Degrees of freedom },
				ylabel={$\| \widetilde{\bm{P}} - \bm{P}\|_{\Le}$},
				xmin=100, xmax=1e5,
				ymin=1e-2, ymax=100,
				xtick={100,1000,10000,1e5},
				ytick={1e-2,1e-1,1, 10,100},
				legend pos= south west,
				ymajorgrids=true,
				grid style=dotted,
				]
				\addplot[
				color=black,
				mark=triangle,
				]
				coordinates {
					( 243 , 1.6492422502470698 )
					( 1125 , 1.1141756364192905 )
					( 3087 , 0.8414250024711007 )
					( 6561 , 0.7180789759606185 )
					( 11979 , 0.6256419808314977 )
				};
				\addlegendentry{$[\Hone{}]^{3\times3}$}
				\addplot[
				color=teal,
				mark=o,
				]
				coordinates {
					( 270 , 1.365295410721133 )
					( 1620 , 0.9348310379951252 )
					( 4914 , 0.7562940652094656 )
					( 11016 , 0.648086737718371 )
					( 20790 , 0.5766243901635305 )
				};
				\addlegendentry{$[\Hc{}]^3$}
				\addplot[
				color=blue,
				mark=square,
				]
				coordinates {
					( 376 , 7.347941949583744e-16 )
				};
				\addlegendentry{$\Hsc{}$}
				\addplot[dashed,color=black, mark=none]
				coordinates {
					(500, 10)
					(25000, 5.2100073095869135)
				};
			\end{loglogaxis}
			\draw (3.,4.7) node[anchor=north west]{$\mathcal{O}(\sqrt{h})$};
		\end{tikzpicture}
	\end{subfigure}
	\begin{subfigure}{0.48\linewidth}
		\begin{tikzpicture}
			\begin{loglogaxis}[
				/pgf/number format/1000 sep={},
				axis lines = left,
				xlabel={Degrees of freedom },
				ylabel={$\|\widetilde{\bm{P}} - \bm{P}\|_{\Hsc{}}$ },
				xmin=100, xmax=1e5,
				ymin=1e-2, ymax=100,
				xtick={100,1000,10000,1e5},
				ytick={1e-2,1e-1,1, 10,100},
				legend pos= south west,
				ymajorgrids=true,
				grid style=dotted,
				]
				\addplot[
				color=black,
				mark=triangle,
				]
				coordinates {
				( 243 , 1.6492422502470698 )
				( 1125 , 1.1141756364192905 )
				( 3087 , 0.8414250024711007 )
				( 6561 , 0.7180789759606185 )
				( 11979 , 0.6256419808314977 )
				};
				\addlegendentry{$[\Hone{}]^{3\times3}$}
				\addplot[
				color=teal,
				mark=o,
				]
				coordinates {
				( 270 , 1.8146543916005167 )
				( 1620 , 1.4630532922020623 )
				( 4914 , 1.3658130115486644 )
				( 11016 , 1.315575081428732 )
				( 20790 , 1.286064678839632 )
				};
				\addlegendentry{$[\Hc{}]^3$}
				\addplot[
				color=blue,
				mark=square,
				]
				coordinates {
				( 376 , 7.357033354028338e-16 )
				};
				\addlegendentry{$\Hsc{}$}
				\addplot[dashed,color=black, mark=none]
				coordinates {
					(500, 0.3)
					(25000, 0.1563002192876074)
				};
				\addplot[dashed,color=black, mark=none]
				coordinates {
					(500, 8)
					(25000, 8)
				};
			\end{loglogaxis}
			\draw (3., 2.55) node[anchor=north west]{$\mathcal{O}(\sqrt{h})$};
			\draw (3.,4.7) node[anchor=north west]{$\mathcal{O}(0)$};
		\end{tikzpicture}
	\end{subfigure}
	\caption{Convergence rates for the discontinuous identity benchmark.}
	\label{fig:dis}
\end{figure}

\begin{figure}
	\centering
	\begin{subfigure}{0.3\linewidth}
		\includegraphics[width=1.0\linewidth]{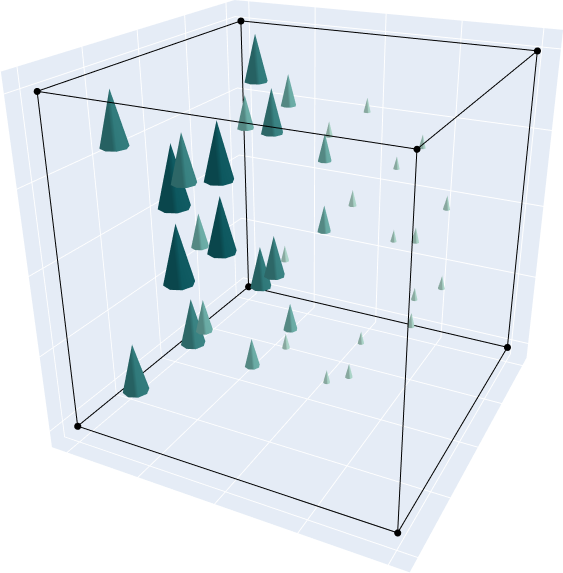}
		\caption{40 elements} 
	\end{subfigure}
	\begin{subfigure}{0.3\linewidth}
		\includegraphics[width=1.0\linewidth]{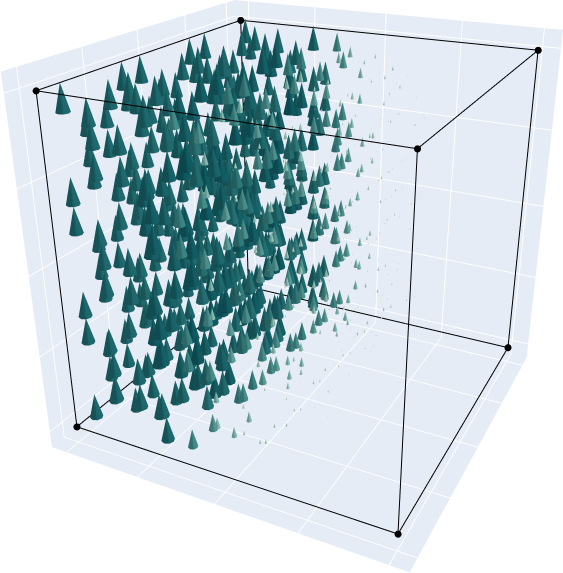}
		\caption{1080 elements}
	\end{subfigure}
	\begin{subfigure}{0.3\linewidth}
		\includegraphics[width=1.0\linewidth]{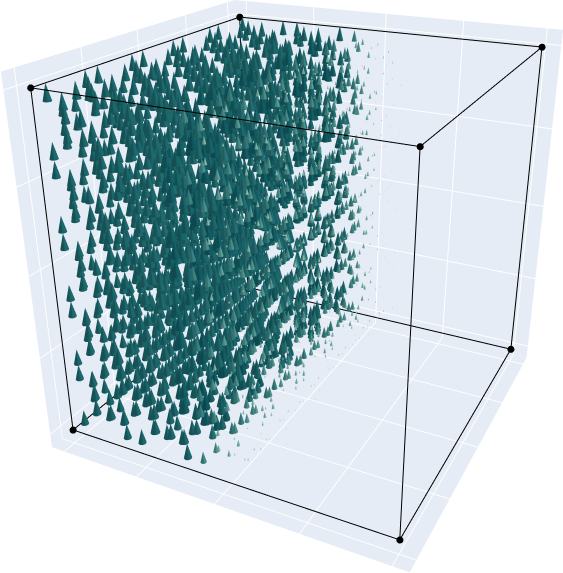}
		\caption{5000 elements}
	\end{subfigure}
	\caption{Discontinuous identity field for various discretizations of the $\Hone$-element displaying the last row of $\Pm$.}
	\label{fig:dispics0}
\end{figure}

\begin{figure}
	\centering
	\begin{subfigure}{0.3\linewidth}
		\includegraphics[width=1.0\linewidth]{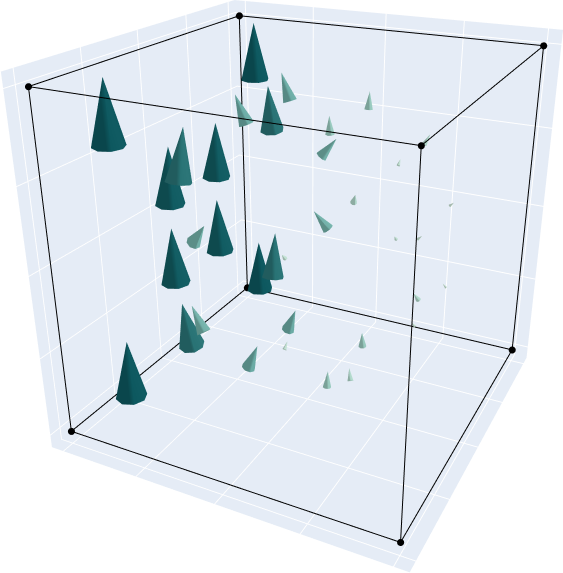}
		\caption{40 elements} 
	\end{subfigure}
	\begin{subfigure}{0.3\linewidth}
		\includegraphics[width=1.0\linewidth]{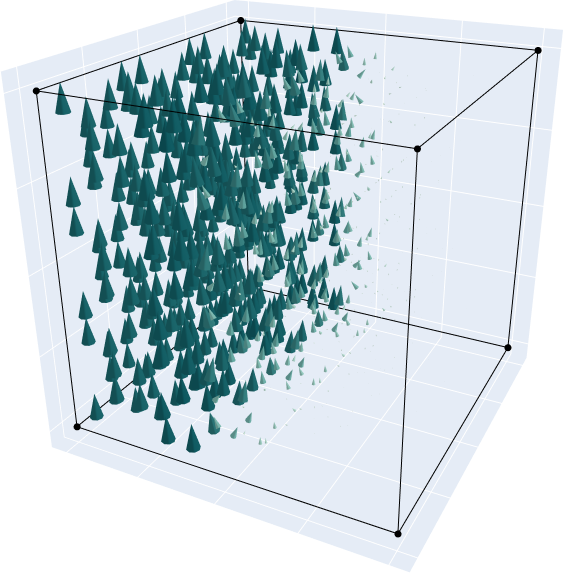}
		\caption{1080 elements}
	\end{subfigure}
	\begin{subfigure}{0.3\linewidth}
		\includegraphics[width=1.0\linewidth]{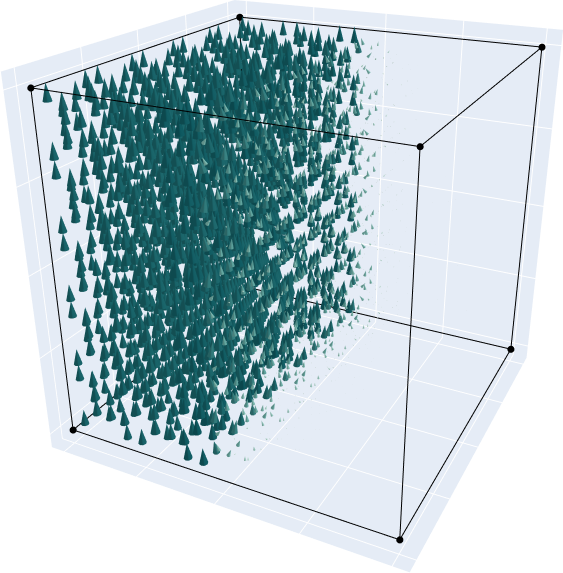}
		\caption{5000 elements}
	\end{subfigure}
	\caption{Discontinuous identity field for various discretizations of the $\Hc{}$-element displaying the last row of $\Pm$.}
	\label{fig:dispics}
\end{figure}

\begin{figure}
	\centering
	\begin{subfigure}{0.3\linewidth}
		\includegraphics[width=1.0\linewidth]{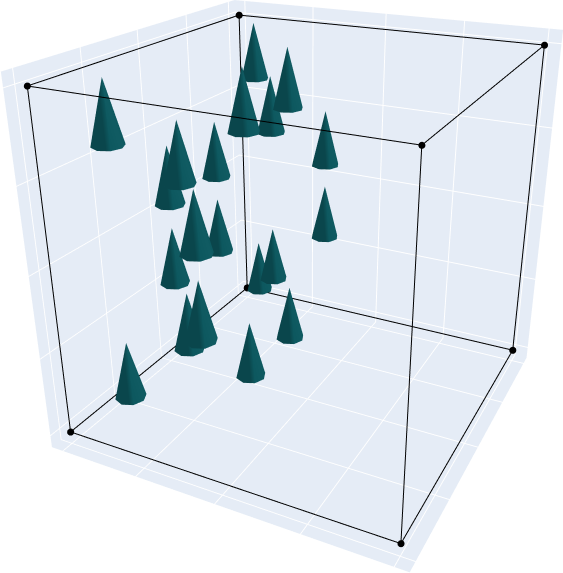}
	\end{subfigure}
	\begin{subfigure}{0.3\linewidth}
		\includegraphics[width=1.0\linewidth]{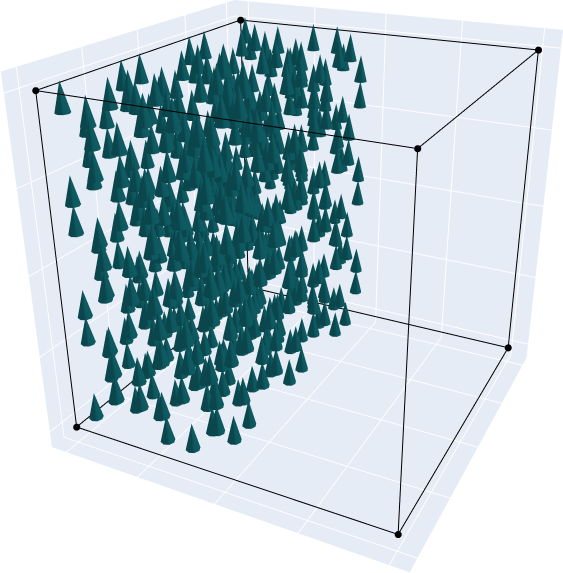}
	\end{subfigure}
	\begin{subfigure}{0.3\linewidth}
		\includegraphics[width=1.0\linewidth]{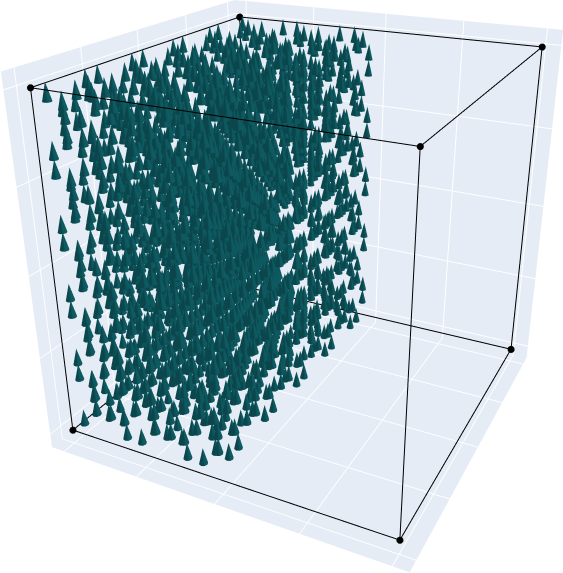}
	\end{subfigure}
	\caption{Discontinuous identity field for various discretizations of the $\Hsc{}$-element displaying the last row of $\Pm$.}
	\label{fig:dispics2}
\end{figure}

As shown in \cref{fig:dis}, the $\Hsc{}$-formulation finds the analytical solution immediately for all domain discretizations, whereas the $\Hone$ and $\Hc{}$ formulations exhibit sub-optimal square root convergence. However, in the $\Hsc{}$-norm only the $\Hone$-formulation continues to converge, while the slope of the $\Hc{}$-formulation quickly tends to zero.
The errors in the solution are clearly visible in the form of noise in \cref{fig:dispics,fig:dispics0}, whereas \cref{fig:dispics2} depicts the discontinuous field as captured by the $\Hsc{}$ formulation. 

\section{Conclusions and outlook}
The relaxed micromorphic model with a symmetric micro-dislocation $\sym \Curl \Pm$ as curvature measure further reduces the continuity requirements of the microdistortion field. As derived by the kernel of the trace and demonstrated by our examples, discontinuous spherical tensors can be captured by the $\Hsc{}$-space, as opposed to the $[\Hone]^{3\times 3}$ and $[\Hc{}]^3$ spaces. In addition, the tests show the corresponding finite element converges optimally for discontinuous spherical fields in $\Hsc{}$, whereas the Lagrangian and N\'ed\'elec elements exhibit sub-optimal square root convergence.
Further, the convergence slope of the N\'ed\'elec-element rapidly flattens in the $\Hsc{}$-norm for discontinuous spherical fields.

These findings serve as a basis for understanding the behaviour of $\Hsc{}$-elements in implementations of the relaxed micromorphic continuum or in computations of the biharmonic equation \cite{Pauly2020}. The relaxed micromorphic complex is needed for future works, where mixed formulations of the relaxed micromorphic model with $\Hsc{}$- and $\HdD{}$-conforming elements are employed to stabilise evaluation with $\Lc \to \infty$. The latter requires the introduction of the $\HdD{}$- or optimally $\HdD{} \cap \Sym(3)$-conforming finite elements.

\section*{Acknowledgements}
We thank Oliver Sander (TU Dresden) for constructive discussions during drafting this manuscript.

Patrizio Neff acknowledges support in the framework of the DFG-Priority Programme
2256 “Variational Methods for Predicting Complex Phenomena in Engineering Structures and Materials”, Neff 902/10-1, Project-No. 440935806. 

\bibliographystyle{spmpsci}   

\setstretch{0.7}
\bibliography{Ref}   

\setstretch{1.}
\appendix

\section{Some mathematical identities}
\label{ap:A}

\subsection{The range of $\sym\Curl$}
In index-notation one finds $\Curl \Pm = - P_{ij,k} \, \varepsilon_{jkl} \, \vb{e}_i \otimes \vb{e}_l$ and as such
\begin{align}
    \sym \Curl \Pm = -\dfrac{1}{2}(P_{ij,k} \, \varepsilon_{jkl} \, \vb{e}_i \otimes \vb{e}_l + P_{ij,k} \, \varepsilon_{jkl} \, \vb{e}_l \otimes \vb{e}_i) \, .
\end{align}
Eliminating the constant and applying the $\Di$ operator yields 
\begin{align}
    -2 \Di \sym \Curl \Pm = (P_{ij,kr} \, \varepsilon_{jkr} \, \vb{e}_i + P_{ij,ki} \, \varepsilon_{jkl} \, \vb{e}_l) = P_{ij,ki} \, \varepsilon_{jkl} \, \vb{e}_l \, ,
\end{align}
since the double contraction between the symmetric and anti-symmetric tensors $P_{ij,kr} \, \varepsilon_{jkr}$ is zero.
The next divergence operator kills the remaining term
\begin{align}
    -2 \di \Di \sym \Curl \Pm = P_{ij,kil} \, \varepsilon_{jkl} = P_{ij,ikl} \, \varepsilon_{jkl} = 0 \, ,
\end{align}
where the change in the order of the partial derivatives is due to Schwarz's theorem.
Therefore, there holds 
\begin{align}
    \di \Di \sym \Curl \Pm = \dfrac{1}{2} \di \Di \Curl \Pm + \dfrac{1}{2} \di \Di (\Curl \Pm)^T = \dfrac{1}{2} \di \Di (\Curl \Pm)^T = 0 \, ,
\end{align}
and 
\begin{align}
    \range(\sym \Curl) \subset \ker(\di \Di) \, .
\end{align}

\subsection{The kernel of $\sym \Curl$}
     The identity $\ker(\Curl) = \range(\D)$ is derived directly by the row-wise application of $\ker(\curl) = \range(\nabla)$.
     For the remaining part we consider $\ker(\sym) = \so$ and as such
     \begin{align}
         \Curl \Pm = \bm{A} \, , \qquad \bm{A} \in \so \, .
     \end{align}
     Further, we can always write $\bm{A} = \Anti(\vb{a})$ for some $\vb{a} \in \mathbb{R}^3$
     and consequently
     \begin{align}
         \Curl \Pm = \Anti(\vb{a}) \, .
         \label{eq:curlanti}
     \end{align}
     Applying the divergence operator on both sides yields
     \begin{align}
         \Di \Anti(\vb{a}) = \Anti(\vb{a})_{,i} \, \vb{e}_i = \vb{a}_{,i} \times \vb{e}_i = - \curl \vb{a} = 0 \, ,
     \end{align}
     since $\Di \Curl \Pm = 0$.
     The latter is equivalent to
     $\vb{a} = \nabla \lambda$ for some scalar field $\lambda :\Omega \mapsto \mathbb{R}$ on the contractible domain $\Omega$.
     Now observe that
     $\Curl(\lambda \one) = -\Anti(\nabla \lambda)$ and thus $\Pm = \lambda \one$ satisfies \cref{eq:curlanti}.
     Clearly, any other field $\bm{T} = \Curl \bm{P}$ where $\bm{T} \notin \so$ is not in $\ker(\sym\Curl)$ simply because $\bm{T}$ is not purely anti-symmetric and consequently, not in $\ker(\sym)$. The remaining part of the kernel is given by gradient fields, finally yielding
     \begin{align}
         \ker(\sym \Curl) = \range(\D) \cup \Sph \, .
     \end{align}

\end{document}